\documentclass[aap,dvips]{arximspdf}
\usepackage{mathbh}
\usepackage{graphicx}

\doi{10.1214/07-AAP465}
\volume{18}
\issue{3}
\pubyear{2008}
\firstpage{967}
\lastpage{996}

\makeatletter

\newtheorem{Thm}{Theorem}[section]
\newtheorem{Lem}[Thm]{Lemma}
\newtheorem{Pte}[Thm]{Proposition}
\newtheorem{Cor}[Thm]{Corollary}
\newproclaim{ex}{Example}
\newproclaim{Rque}{Remark}
\makeatother

\begin{document}
\begin{frontmatter}

\title{Proliferating parasites in dividing cells: Kimmel's branching
model revisited}
\runtitle{Proliferating parasites in dividing cells}

\begin{aug}
\author[A]{\fnms{Vincent} \snm{Bansaye}\ead[label=e1]{bansaye@ccr.jussieu.fr}\corref{}}
\runauthor{V. Bansaye}
\affiliation{Universit\'e Paris 6}
\address[A]{UPMC et C.N.R.S. UMR 7599\\
175 rue du Chevaleret \\
75 013 Paris\\
France\\
\printead{e1}} 
\end{aug}

\received{\smonth{2} \syear{2007}}
\revised{\smonth{2} \syear{2007}}

%
\begin{abstract}
We consider a branching model introduced by Kimmel for cell division
with parasite infection. Cells
contain proliferating parasites which are shared randomly between the
two daughter cells when they divide.
We determine the probability that the organism recovers, meaning that
the asymptotic proportion
of contaminated cells vanishes. We study
the tree of contaminated cells, give the asymptotic number of
contaminated cells and the asymptotic
proportions of contaminated cells with a given number of parasites. This
depends on domains inherited from the behavior of branching processes
in random environment (BPRE) and given
by the bivariate value of the means of parasite offsprings. In one of
these domains, the convergence of proportions
holds in probability, the limit is deterministic and given by the
Yaglom quasistationary
distribution. Moreover, we get an interpretation of the limit of the Q-process
as the size-biased quasistationary distribution.
\end{abstract}

%
\pdfkeywords{60J80, 60J85, 60K37, 92C37, 92D25, 92D30,
Bienayme Galton Watson process (BGW), branching processes in random environment (BPRE),
Markov chain indexed by a tree, quasistationary distribution, empirical measures,}

\end{frontmatter}

\section{Introduction}

We consider the following model for cell division with parasite
infection. Unless
otherwise specified, we start with a single cell infected with a single
parasite. At
each generation, each parasite multiplies independently, each cell
divides into two daughter cells and
the offspring of each parasite is shared independently into the two
daughter cells.
It is convenient to distinguish a first daughter cell called $0$ and a
second one called $1$
and to write $Z^{(0)}+Z^{(1)}$ the number of offspring of a parasite,
$Z^{(0)}$ of which go
into the first daughter cell and $Z^{(1)}$ of which into the second
one. The symmetric sharing is the
case when $(Z^{(0)},Z^{(1)})\stackrel{d}{=}(Z^{(1)},Z^{(0)})$.
Even in that case, the
sharing of parasites can be unequal [e.g., when
$\mathbb{P}(Z^{(0)}Z^{(1)}=0)=1$].

\begin{figure}

\includegraphics{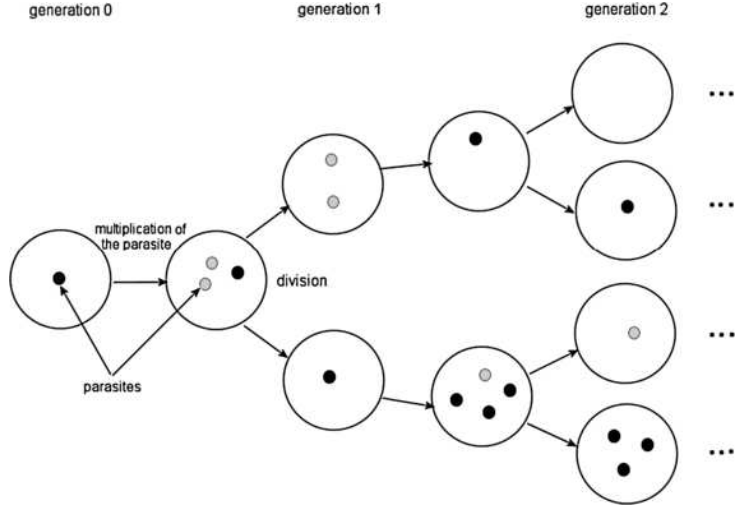}

 \caption{Multiplication of parasites and cell division.\label{f1}}
\end{figure}

We denote by $\mathbb{T}$ the binary genealogical tree of the cell
population, by $\mathbb{G}_n$
(resp. $\mathbb{G}_n^*$) the set of cells at generation
$n$ (resp. the set of contaminated cells at generation~$n$) and by
$Z_{\mathbf{i}}$ the number of parasites of cell
$\mathbf{i}\in\mathbb{T}$, that is,
\[
\mathbb{G}_n:=\{0,1\}^{n}, \qquad
\mathbb{G}_n^*:=\{\mathbf{i} \in\mathbb{G}_n \dvtx Z_{\mathbf{i}}>0\}, \qquad
\mathbb{T}:=\bigcup_{n\in\mathbb{N}} \mathbb{G}_n.
\]
For every cell $\mathbf{i}\in\mathbb{T}$,
conditionally on $Z_{\mathbf{i}}=x$, the numbers of parasites
$(Z_{\mathbf{i}0}, Z_{\mathbf{i}1})$ of its two daughter cells is
given by
\[
\sum_{k=1}^x \bigl(Z^{(0)}_k (\mathbf{i}), Z^{(1)}_k(\mathbf{i})\bigr),
\]
where $(Z_k^{(0)} (\mathbf{i}), Z_k^{(1)} (\mathbf{i}))_{\mathbf
{i} \in\mathbb{T}, k\geq1}$ is an i.i.d. sequence
distributed as $(Z^{(0)}, Z^{(1)})$ (see Figure~\ref{f1}).

This is a discrete version of the model introduced by
Kimmel in \cite{kim}. In particular, it contains the following model
with binomial repartition of parasites. Let $Z$ be a random
variable in $\mathbb{N}$ and $p\in[0,1]$. At each generation, every
parasite
multiplies independently with the same reproduction law $Z$. When the
cells divides,
every parasite chooses independently the first daughter cell
with probability
$p$ (and the second one with probability $1-p$). It contains also the
case when every parasite
gives birth to a random cluster of parasites of size $Z$ which goes to
the first cell with
probability $p$ (and to the second one with probability $1-p$).

We introduce for $a\in\{0,1\}$
%
\begin{equation}\label{nota}
m_a := \mathbb{E}\bigl(Z^{(a)}\bigr) \qquad\forall s\geq0,
f_a(s):=\mathbb{E}\bigl(s^{Z^{(a)}}\bigr).
\end{equation}
We assume $0<m_0<\infty$,
$0<m_1<\infty$ and to avoid trivial cases, we require
%
\begin{equation}\label{condp}
\mathbb{P}\bigl(\bigl(Z^{(0)}, Z^{(1)}\bigr)=(1,1)\bigr) <1, \qquad
\mathbb{P}\bigl(\bigl(Z^{(0)},Z^{(1)}\bigr) \in\{(1,0),(0,1)\}\bigr)< 1.
\end{equation}

This model is a Markov chain indexed by a tree. This subject has
been studied in the literature (see e.g., \cite{atarbre,benj})
in the symmetric independent case. In this case, for every
$(\mathbf{i},k) \in\mathbb{T}\times\mathbb{N}$, we have
\[
\mathbb{P}\bigl((Z_{\mathbf{i}0}, Z_{\mathbf{i}1})=(k_0,k_1)\mid Z_{\mathbf{i}}=k\bigr)
=\mathbb{P}(Z_{\mathbf{i}0}=k_0\mid Z_{\mathbf{i}}=k)
\mathbb{P}(Z_{\mathbf{i}0}=k_1\mid Z_{\mathbf{i}}=k)
\]
which require
that $Z^{(0)}$ and $Z^{(1)}$ are i.i.d. in this model. Guyon
\cite{guy} studies a Markov chain indexed by a binary tree where
asymmetry and dependence are allowed and limit theorems are
proved. But the case where his results apply is degenerate (this
is the case $m_0m_1\leq1$ and the limit of the number of
parasites in a random cell line is zero). Moreover, adapting his
arguments for the theorems stated here appears to be cumbersome
(see the remark in Section \ref{s5.2} for details). In the same vein, we
refer to
\cite{ewans,tad} (cellular aging).

The total population of parasites at generation $n$, which we denote
by $\mathcal{Z}_n$, is
a Bienaym\'{e} Galton--Watson process (BGW) with
reproduction law $Z^{(0)}+Z^{(1)}$. We call Ext (resp.
$\operatorname{Ext}^c$) the event extinction of the parasites (resp.
nonextinction of the parasites),
%
\begin{eqnarray}\label{defext}
\mathcal{Z}_n &=& \sum_{\mathbf{i}\in\mathbb{G}_n} Z_{\mathbf{i}},
\nonumber\\
\mathrm{Ext} &=& \{\exists n \in\mathbb{N}\dvtx \mathcal{Z}_n=0\},
\\
\mathrm{Ext}^c&=& \{\forall n\in\mathbb{N}\dvtx \mathcal{Z}_n>0\}.
\nonumber
\end{eqnarray}
Another process that appears naturally is the number
of parasites in a random cell line. More precisely, let $(a_i)_{i
\in\mathbb{N}}$ be an i.i.d. sequence independent of
$(Z_{\mathbf{i}})_{\mathbf{i}\in\mathbb{T}}$ such that
%
\begin{equation}\label{iid}
\mathbb{P}(a_1=0)=\mathbb{P}(a_1=1)=1/2.
\end{equation}
Then $(Z_n)_{n\in\mathbb{N}}=(Z_{(a_1,a_2,\ldots,a_n)})_{n \in\mathbb{N}}$ is a
Branching Process in Random Environment (BPRE).

The first question we answer here arose from observations made by
de Paepe, Paul and Taddei at TaMaRa's Laboratory
(H\^{o}pital Necker, Paris). They have infected the bacteria \textit{E.
coli} with a parasite (lysogen bacteriophage M13). A fluorescent
marker allows them to see the level of contamination of cells.
They observed that a very contaminated cell often gives birth to a
very contaminated cell which dies fast and to a much less
contaminated cell whose descendance may survive. So cells tend to
share their parasites unequally when they divide so that there are
lots of healthy cells. This is a little surprising since one
could think that cells share equally all their biological content
(including parasites). In Section~\ref{s3}, we prove that if $m_0m_1\leq
1$, the organism recovers a.s. (meaning that the number of
infected cells becomes negligible compared to the number of cells
when $n \rightarrow\infty$). Otherwise
the organism recovers iff parasites die out (and the probability is
less than~$1$).
%

In Section \ref{s4}, we consider the tree of contaminated cells. We denote by
$\partial\mathbb{T}$ the boundary of the cell tree $\mathbb{T}$ and
by $\partial \mathbb{T}^*$ the infinite
lines of contaminated cells, that~is
\[
\partial\mathbb{T}=\{0,1\}^{\mathbb{N}}, \qquad
\partial\mathbb{T}^*=\{ \mathbf{i} \in
\partial\mathbb{T}\dvtx \forall n \in\mathbb{N}, Z_{\mathbf{i}\vert n}\ne0\}.
\]
We shall prove that
the contaminated cells are not concentrated in a cell line. Note that
if $m_0+m_1>1$, conditionally on $\operatorname{Ext}^c$,
$\partial\mathbb{T}^*\ne\varnothing$ since at each generation, one
can choose a daughter cell whose parasite descendance
does not become extinct.

\begin{figure}

\includegraphics{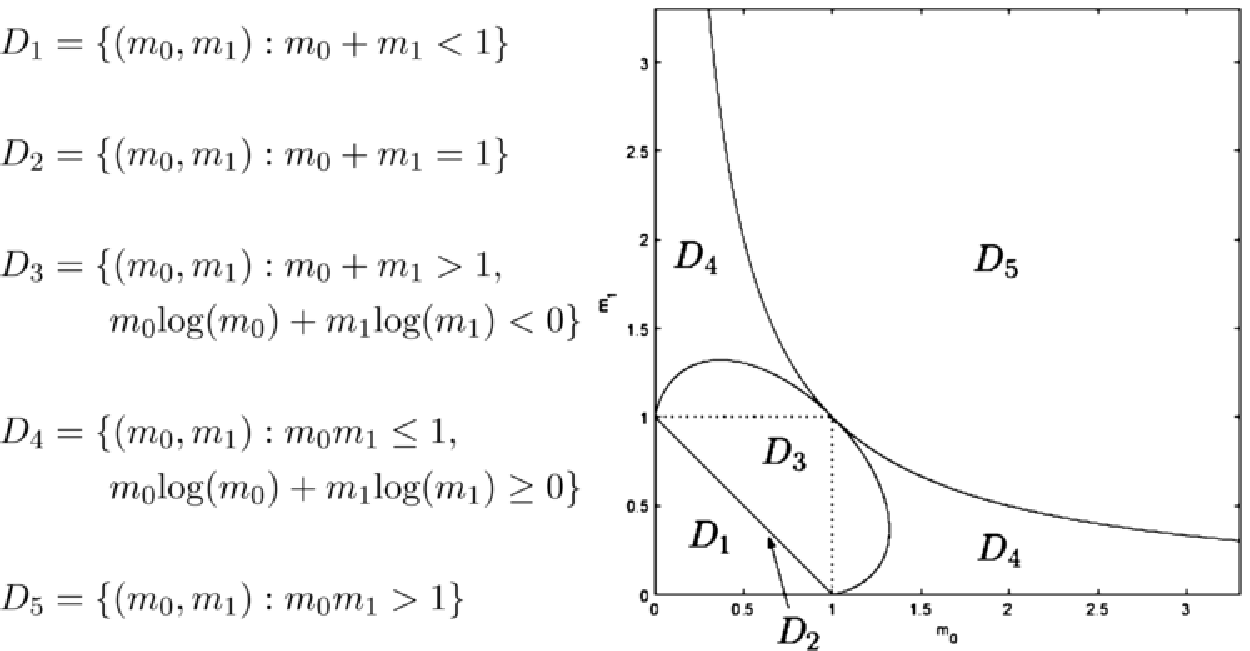}

  \caption{\label{f2}}
\end{figure}

The rest of the work is devoted to the convergence of the number
of contaminated cells in generation $n$ and the convergence of
proportions of contaminated cells with a given number of
parasites (Section~\ref{s5}). These asymptotics depend on $(m_0,m_1)$
and we distinguish five different cases which come from the
behavior of the BGW process $\mathcal{Z}_n$ and the BPRE $Z_n$
(Section~\ref{s2}), shown in Figure~\ref{f2}.

If $(m_0,m_1)\in D_5$, the contaminated cells become largely
infected (Theorem \ref{limf3}). The main two results correspond to
cases $(m_0,m_1)\in D_3$ and $(m_0,m_1)\in D_1$ and are given by the
following two theorems.
\begin{Thm}\label{thm5.2}
If $(m_0,m_1) \in D_3$, conditionally on $\operatorname{Ext}^c$,
the following convergence holds in probability for every $k\in\mathbb{N}$,
\[
\# \{ \mathbf{i}\in\mathbb{G}_n ^*\dvtx
Z_{\mathbf{i}}=k\}/\# \mathbb{G}_n^*
\stackrel{n \rightarrow\infty}{\longrightarrow} \mathbb{P}(\Upsilon =k),
\]
\end{Thm}
where $\Upsilon $ is the Yaglom
quasistationary distribution of the BPRE $(Z_n)_{n\in\mathbb{N}}$ (see
\cite{AN,bpree}). Note that the limit is deterministic and depends
solely on the marginal laws of $(Z^{(0)},Z^{(1)})$ (see
Proposition \ref{yag}). This gives then a way to compute $\Upsilon $ as
a deterministic limit, although it is defined by conditioning on a
vanishing event. Kimmel \cite{kim} considers the symmetric case
($(Z^{(0)},Z^{(1)})\stackrel{{d}}{=}(Z^{(1)},Z^{(0)})$)
with $m_0=m_1< 1 <m_0+m_1$ in a continuous analogue of this model
(cells divide after an exponential time). The
counterpart of his result in the discrete case is easy to prove
[see (\ref{kim})] and makes a first link with $\Upsilon $.
\[
\lim_{n\rightarrow\infty}\mathbb{E}(\# \{ \mathbf{i}\in\mathbb{G}_n \dvtx
Z_{\mathbf{i}}=k\})/\mathbb{E}(\# \mathbb{G}_n^*)= \mathbb{P}(\Upsilon=k).
\]
Moreover, the proportions of contaminated cells on the boundary
of the tree whose ancestors at generation $n$ have a given number of parasites
converge to the size-biased distribution of $\Upsilon $ letting
$n\rightarrow\infty$ (Corollary \ref{limf13}). This gives a
pathwise interpretation that the limit of the Q-process
associated to $Z_n$ (see \cite{afa2,AN}) is the size-biased
quasistationary distribution.
\begin{Thm}\label{thm5.7}
If $(m_0,m_1) \in D_1$,
$(\#\{\mathbf{i} \in\mathbb{G}_n^* \dvtx Z_{\mathbf{i}}=k \})_{k \in\mathbb{N}}$
conditioned on $\mathcal{Z}_n>0$ converges in distribution as
$n\rightarrow \infty$
to a finite random sequence ($N_k)_{k \in\mathbb{N}}$.
\end{Thm}

We obtain a similar result in the case $(m_0,m_1)\in D_2$
(Theorem \ref{limf4}) and we get the following
asymptotics (Theorem \ref{Text} and Corollaries \ref{gn}, \ref{ggn},
\ref{gggn}).

If $(m_0,m_1)\in D_3$ (resp. $D_5$), then
conditionally on $\operatorname{Ext}^c$, $\# \mathbb{G}^*_n/(m_0+m_1)^n$
(resp. $\#\mathbb{G}_n^*/2^n$) converges in probability to a finite
positive r.v.

If $(m_0,m_1)\in D_1$ (resp. $D_2$), then
$\# \mathbb{G}^*_n$ (resp. $\# \mathbb{G}_n^*/n$)
conditioned by $\# \mathbb{G}_n^*>0$ converges in distribution to a finite
positive r.v.

In the case $(m_0,m_1) \in D_4$, we get only some estimates of the
asymptotic of $\# \mathbb{G}^*_n$ which are different from those which
hold in the other domains. Our conjecture is that $\# \mathbb{G}^*_n$
has also a deterministic asymptotic, which depends on three
subdomains (the interior of $D_4$ and its boundaries). As a
perspective, we are also interested in determining which types of
convergences hold in $D_4$ for the proportions of contaminated
cells with a given number of parasites (see Section~\ref{s5.5}).

Moreover, we wonder if the convergences stated above hold a.s. and if
they extend to
the continuous case and complement the results of Kimmel. Finally, in a
work in progress with Beresticky and Lambert,
we aim at determining the localizations of contaminated cells and the
presence of
cells filled-in by parasites on the boundary of the tree (branching
measure and multifractal analysis).

\section{Preliminaries}\label{s2}

In this section, we give some useful results about the two
processes introduced above. First define:
%
\begin{equation}\label{lettre}
m:=\tfrac{1}{2} (m_0+m_1).
\end{equation}
We use the classical notation, where
for every $\mathbf{i}=(\alpha_1,\ldots,\alpha_{n}) \in\mathbb{G}_n$,
\begin{eqnarray*}
\vert\mathbf{i} \vert &=& n, \qquad
\mathbf{i} \vert k=(\alpha_1,\ldots,\alpha_{k})
\qquad \mbox{for every } k\leq n, \qquad
\\
\mathbf{j} &<& \mathbf{i} \qquad \mbox{if } \exists k<n \dvtx \mathbf{i}\vert k= \mathbf{j}.
\end{eqnarray*}

\subsection{Results on the BGW process
$(\mathcal{Z}_n)_{n\in\mathbb{N}}$}\label{s2.1}

The results stated hereafter are well known and can be found in
\cite{AN}. First, the probability
of extinction of the parasites satisfies
\begin{eqnarray*}
\mathbb{P}(\mathrm{Ext}) &=& \operatorname{inf}\bigl\{s \in[0,1] \dvtx
\mathbb{E} \bigl(s^{Z^{(0)}+Z^{(1)}}\bigr)=s\bigr\};
\\
\mathbb{P}(\operatorname{Ext}) &=& 1 \qquad \mbox{iff } m_0+m_1\leq1/2.
\end{eqnarray*}

From now, we assume
\[
\check{m}:=\mathbb{E}\bigl(\bigl(Z^{(0)}+Z^{(1)}\bigr)
\log ^+\bigl(Z^{(0)}+Z^{(1)}\bigr)\bigr) <\infty.
\]
Then there exists a random variable $W$ such that
%
\begin{equation}\label{nbparas}
\frac{\mathcal{Z}_n}{(m_0+m_1)^n}\stackrel{n\rightarrow\infty}{\longrightarrow} W,
\qquad
\mathbb{P}(W=0)=\mathbb{P}(\operatorname{Ext}),
\qquad
\mathbb{E}(W)=1.
\end{equation}

In the case $m_0+m_1<1$, there exists $b>0$ such that
$\mathbb{P}(\mathcal{Z}_n>0)\stackrel{n\rightarrow\infty}{\sim}b
(m_0+m_1)^n$. Then, there exists $U>0$ such that
%
\begin{equation}\label{mintot}
\mathbb{P}(\mathcal{Z}_n>0)\geq U(m_0+m_1)^n.
\end{equation}
Moreover $(\mathcal{Z}_n)_{n\in\mathbb{N}}$ conditioned to be
nonzero converges to
a variable called the Yaglom quasistationary distribution and we set
%
\begin{equation}\label{yagnb}
\mathcal{B}(s):=\lim_{n\rightarrow\infty}
\mathbb{E}(s^{\mathcal{Z}_n} \mid\mathcal{Z}_{n}>0).
\end{equation}
We consider then $\mathcal{B}_{n,k}(s):=\mathbb{E}(s^{\mathcal{Z}_n}
\mid\mathcal{Z}_{n+k}>0)$
which satisfies 
%
\begin{equation}\label{limunq}
\lim_{n\rightarrow\infty} \mathcal{B}_{n,k}(s)
= \frac{\mathcal{B}(s)-\mathcal{B}(sf_k(0))}{1-\mathcal{B}(f_k(0))}.
\end{equation}
Moreover $\mathcal{B}$ is differentiable at $1$ (Lemma $1$ on page 44
in \cite{AN}) and we get
%
\begin{equation}\label{limlim}
\lim_{k\rightarrow\infty} \lim_{n\rightarrow\infty} \mathcal{B}
_{n,k}(s)=\frac{s\mathcal{B}'(s)}{\mathcal{B}'(1)}.
\end{equation}
This is the probability generating function of the size-biased Yaglom
quasistationary distribution, which is also the stationary
distribution of the Q-process.

Finally if $\hat{m}:=\mathbb{E}
((Z^{(0)}+Z^{(1)})((Z^{(0)}+Z^{(1)})-1))<\infty$ and $2m\ne1$, then
%
\begin{equation}\label{mom2GW}
\mathbb{E}\bigl(\mathcal{Z}_n(\mathcal{Z}_n-1)\bigr)
= \hat{m}(2m)^n \frac{(2m)^{n}-1}{(2m)^2-2m}.
\end{equation}
%

\subsection{Properties of the BPRE $(Z_n)_{n\in\mathbb{N}}$}\label{s2.2}

Recall that $(Z_n)_{n\in\mathbb{N}}$ is the population of parasites
in a
uniform random cell line. Then $(Z_n)_{n\in\mathbb{N}}$
is a BPRE with two equiprobable environments. More precisely, for
each $n\in\mathbb{N}$, conditionally on $a_n=a$ with $a\in\{0,1\}$ [see
(\ref{iid})], all parasites behave independently of one another
and each of them gives birth to $Z^{(a)}$ children. The size of
the population at generation~$0$ is denoted by $k$ and we note
$\mathbb{P}_k$ the associated probability. Unless otherwise mentioned, the
initial state is equal to $1$. For the general theory, see, for example,
\cite{Dek,bpree,Guiv,bpre}. In the case
$Z^{(0)}\stackrel{{d}}{=}Z^{(1)}$, $(Z_n)_{n\in\mathbb{N}}$
is a BGW with reproduction law~$Z^{(0)}$.

For $\mathbf{i}=(\alpha_1,\ldots,\alpha_n) \in\mathbb{G}_n$, we define
\[
f_{\mathbf{i}}:=f_{\alpha_1}\circ\cdots\circ f_{\alpha_n}, \qquad
m_{\mathbf{i}}=\prod_{i=1}^{n} m_{\alpha_i},
\]
and for all $(n,k) \in\mathbb{N}\times\mathbb{N}^*$ and $\mathbf
{i} \in\mathbb{G}_n$,
\[
\mathbb{E}_k \bigl(s^{Z_n}\mid(a_1,\ldots,a_{n})
=\mathbf{i}\bigr)=f_{\mathbf{i}}(s)^k.
\]
Then for all $(n,k) \in\mathbb{N}\times\mathbb{N}^*$ and $s\in
[0,1]$,
%
\begin{equation}\label{eqfond}
\mathbb{E}_k (s^{Z_n} )=2^{-n}\sum_{\mathbf{i}\in \mathbb{G}_n}
f_{\mathbf{i}}(s)^k.
\end{equation}

First, for every $n\in\mathbb{N}$, $\mathbb{E}(Z_{n+1}\mid
Z_n)=mZ_n$ and $\mathbb{E}(Z_n)=m^n$.

Moreover, as $(\mathbb{P}(Z_n=0))_{n\in\mathbb{N}}$ is an increasing
sequence, it
converges to the probability of extinction $p$ of the process.
Recalling (\ref{nota}), we have the following result (see
\cite{bpre} or \cite{at}).
\begin{Pte}\label{ext}
If $m_0m_1\leq1$, then $p=1$. Otherwise $p<1$.
\end{Pte}

In the subcritical case ($m_0m_1< 1$), the process $Z_n$
conditioned to be nonzero which is denoted by $Z_n^*$ converges
weakly (Theorem 1.1 in \cite{bpree}). By analogy with BGW, we call
its limit distribution the Yaglom quasistationary distribution and
denote it by $\Upsilon $. That is,
\[
\forall s\in[0,1] \qquad
\mathbb{E} (s^{Z_n}\mid Z_n>0 )
\stackrel{n\rightarrow\infty}{\longrightarrow}
\mathbb{E}(s^{\Upsilon })=G(s).
\]
In the subcritical case, the asymptotics of
$(\mathbb{P}(Z_n>0))_{n\in\mathbb{N}}$ when $n$ is large depends on
the sign of
$m_0\log (m_0)+m_1\log (m_1)$ (see \cite{bpree}).
Now, we require also that
%
\begin{equation}\label{conH}
m_0\log (m_0)+m_1\log (m_1)<0;
\qquad
\mathbb{E}(Z_{a}\log ^+(Z_a))<\infty.
\end{equation}
Then we say that $Z_n$ is
strongly subcritical and there exists $c>0$ such that as $n$
tends to $\infty$ (Theorem 1.1 in \cite{bpree}),
%
\begin{equation}\label{equiv}
\mathbb{P}(Z_n>0)\sim c m^n.
\end{equation}
Moreover, in that case, $\Upsilon $ is characterized by
\begin{Pte}\label{yag}
$G$ is the unique probability generating function which satisfies
%
\begin{eqnarray}\label{equatfonct}
G(0) &=& 0, \qquad
G'(1)<\infty, \qquad
\nonumber\\[-8pt]
\\[-8pt]
\frac{G(f_0(s))+G(f_1(s))}{2} &=& mG(s)+(1-m).
\nonumber
\end{eqnarray}
\end{Pte}


To prove the uniqueness, we need and prove below the following result.
\begin{Lem}\label{lemfonct}
If $H\dvtx [0,1]\mapsto\mathbb{R}$ is continuous, $H(1)=0$
and
%
\begin{equation}\label{equatder}
H=\frac{H\circ f_0\cdot f_0'+H\circ f_1 \cdot f_1'}{2m},
\end{equation}
then $H=0$.
\end{Lem}
\begin{pf*}{Proof of Proposition \ref{yag}}
The finiteness of $G'(1)=\mathbb{E}(\Upsilon )$ is the second part of
Theorem 1.1 in \cite{bpree}.

We characterize now the probability generating function $G$ of the limit
distribution:
\begin{eqnarray*}
&& 1-\mathbb{E} (s^{Z_{n+1}}\mid Z_{n+1}>0)
\\
&&\qquad = \frac{1-\mathbb{E} (s^{Z_{n+1}} )}{\mathbb{P}(Z_{n+1}>0)}
=\frac{1}{\mathbb{P}(Z_{n+1}>0)}
\sum_{i=1}^{\infty} \mathbb{P}(Z_{n}=i)
\bigl(1-\mathbb{E}_i (s^{Z_1} ) \bigr)
\\
&&\qquad = \frac{\mathbb{P}(Z_{n}>0)}{\mathbb{P}(Z_{n+1}>0)}
\frac{1}{\mathbb{P}(Z_{n}>0)} \sum_{i=1}^{\infty}
\mathbb{P}(Z_{n}=i) \biggl(1-\frac{f_0(s)^i+f_1(s)^i}{2} \biggr)
\\
&&\qquad =\frac{\mathbb{P}(Z_{n}>0)}{\mathbb{P}(Z_{n+1}>0)}
\bigl(1-\mathbb{E}\bigl(f_0(s)^{Z_{n}}\mid Z_{n}>0\bigr)
+ 1-\mathbb{E} \bigl(f_1(s)^{Z_{n}}\mid Z_{n}>0\bigr) \bigr)/2.
\end{eqnarray*}
And (\ref{equiv}) ensures that $\mathbb{P}(Z_{n}>0)/\mathbb{P}(Z_{n+1}>0) \stackrel
{ \scriptscriptstyle n\rightarrow\infty} {\longrightarrow}m^{-1}$,
so that
\[
1-G(s)=\frac{1-G(f_0(s))  +  1-G(f_1(s))}{2m}.
\]

Finally we prove the uniqueness of solutions of this equation. Let
$G$ and $F$ be two probability generating functions which are
solutions of $(\ref{equatfonct})$. Choose \mbox{$\alpha>0$} such that
$G'(1)=\alpha F'(1)$. Putting $H:=G-\alpha F$, $H'$ satisfies
equation $(\ref{equatder})$ and $H'(1)=0$. Thus Lemma
$\ref{lemfonct}$ gives $H'=0$. As $H(0)=0$, $H=0$. Moreover,
$F(1)=G(1)=1$, so $\alpha=1$ and $F=G$.
\end{pf*}
\begin{pf*}{Proof of Lemma \ref{lemfonct}}
If $H\ne0$ then there exists $\alpha\in[0,1[$
such that
\[
\beta:=\sup \{| H(s)|\dvtx s \in[0, \alpha]\}\ne0.
\]
Let $\alpha_n \in[0,1[$ such that $\alpha_n\stackrel{n\rightarrow
\infty}{\longrightarrow}1$ and $\alpha\leq\alpha_n \leq1$.
Then, for every $n\in\mathbb{N}$, there exists $\beta_n \in
[0,\alpha_n]$ such that:
\begin{eqnarray*}
\sup \{| H(s)|\dvtx s \in[0, \alpha_n]\}&=& \vert H(\beta_n)\vert
\\
&\leq& \frac{\vert H(f_0(\beta_n)) \vert f_0'(\beta_n)
+ \vert H(f_1(\beta_n)) \vert f_1'(\beta_n)}{2m}
\\
&<& \sup \{|H(s)|\dvtx s \in[0,1]\},
\end{eqnarray*}
since
$\sup \{| H(s)|\dvtx s \in[0, 1]\}\ne0$ and
$(2m)^{-1}(f_0'(\beta_n)+f_1'(\beta_n))<1$. As $I\cap
J=\varnothing$, $\sup I <\sup(I\cup J) \Rightarrow\sup I<\sup J$, we get
\[
\beta\leq\sup \{| H(s)|\dvtx s \in[0, \alpha_n]\}
< \sup \{| H(s)| \dvtx s \in\,]\alpha_n, 1]\}.
\]
And $H(s) \stackrel{s\rightarrow1}{\longrightarrow}0$ leads to a
contradiction letting $n\rightarrow\infty$. So $H=0$.
\end{pf*}
%

In the subcritical case ($m_0m_1\leq1$), if
$m_0\log (m_0)+m_1\log (m_1)>0$ [resp.
$m_0\log (m_0)+m_1\log (m_1)=0$], we say that $Z_n$
is weakly
subcritical (resp. intermediate subcritical) and we have
$\mathbb{P}(Z_n>0)\sim c' n^{-3/2}\gamma^n$ [resp. $\mathbb{P}(Z_n>0)\sim
c''n^{-1/2}m^n$] for some $\gamma<m, c'>0,c''>0$ (see \cite{bpree} for details).

Finally we have the following expected result in the supercritical
case \cite{atr}.
\begin{Pte}\label{explo}
If $m_0m_1> 1$, $\mathbb{P}(Z_n\stackrel{n\rightarrow
\infty}{\longrightarrow}\infty \mid \forall n \in\mathbb{N}\dvtx Z_n>0 )=1$.
\end{Pte}
%

\section{Probability of recovery}\label{s3}

We say that the organism recovers if the number of
contaminated cells becomes negligible compared to the number of cells
when $n \rightarrow\infty$. We determine here the
probability of this event. Actually if this probability is not
equal to $1$, then
the parasites must die out for the organism to recover.
\begin{Thm}\label{Text}
There exists a random variable $L\in[0,1]$ such that
\[
\#\mathbb{G}_n^*/2^n \stackrel{n\rightarrow\infty}{\longrightarrow}L.
\]
If $m_0m_1\leq1$ then $\mathbb{P}(L=0)=1$.

Otherwise $\mathbb{P}(L=0)<1$ and $\{L=0\}=\operatorname{Ext}$.
\end{Thm}
\begin{Rque*}%
In the case $m_0+m_1>1$ and $m_0m_1\leq1$, the population of
parasites may explode although the organism recovers.

This theorem states how unequal the sharing of parasites must be
for the organism to recover. More precisely, let $m_0=\alpha M, \
m_1=(1-\alpha)M$ where $M>0$ is the parasite growth rate. Then the
organism recovers a.s. iff
\[
M\leq2 \quad \mbox{or} \quad
\alpha\notin \bigl]\bigl(1-\sqrt{1-4/M^2}\bigr)/2,
\bigl(1+\sqrt{1-4/M^2}\bigr)/2\bigr[ \qquad  (M>2).
\]
%
\end{Rque*}

Note that for all $n\in\mathbb{N}$,
\[
\mathbb{E} \biggl(\frac{\# \mathbb{G}^*_n }{2^n} \biggr)
= \frac{\mathbb{E}(\sum_{i\in\mathbb{G}_n}
\mathbh{1}_{Z_{\mathbf{i}}>0})}{2^n}=\mathbb{P}(Z_n>0).
\]
Recalling that $p$ is the probability of extinction of
$(Z_n)_{n\in\mathbb{N}}$,
%
\begin{equation}\label{esp}
\forall n\in\mathbb{N} \qquad
\mathbb{E} \biggl(\frac{\#\mathbb{G}^*_n }{2^n} \biggr)
= \mathbb{P}(Z_n>0)\stackrel{n\rightarrow\infty}{\longrightarrow}1-p.
\end{equation}
The last equality gives also the asymptotic of
$\mathbb{E}(\# \mathbb{G}^*_n )$
as $n\rightarrow\infty$ in the case $m_0m_1<1$ [see Section~\ref{s2.2}
for the asymptotic of $\mathbb{P}(Z_n>0$), which depends on the sign of
$m_0\log (m_0)+m_1\log (m_1)$] and in the case
$m_0m_1=1$ (see \cite{afa,koz}).
%
%
\begin{pf*}{Proof of Theorem \ref{Text}}
As $\#\mathbb{G}_n^*/2^n$ decreases as $n$ increases, it converges as
$n\rightarrow\infty$.

Monotone convergence of $\# \mathbb{G}_n^*/2^n$ to $L$ as $n\rightarrow\infty$
and (\ref{esp}) ensure that $\mathbb{E}(L)=1-p$. Using Proposition
\ref{ext}, we get $\mathbb{P}(L=0)=1$ iff $m_0m_1\leq1$.

Obviously $\{L=0\}\supset\operatorname{Ext}$.
Denote by $\mathcal{P}(n)$ the set of parasites at generation $n$ and
for every
$\mathbf{p} \in\mathcal{P}(n)$, denote by $N_k(\mathbf{p})$ the
number of
cells at generation $n+k$ which contain at least a parasite whose
ancestor is $\mathbf{p}$. Then, for every $n\in\mathbb{N}$,
\[
\{L=0\}=\bigcap_{\mathbf{p} \in\mathcal{P}(n)}
\biggl \{\frac{N_k(\mathbf{p})}{2^k}\stackrel{k\rightarrow\infty}
{\longrightarrow}0 \biggr\}.
\]
As $T_n:=\operatorname{inf}\{k\geq0 \dvtx\mathcal{Z}_k \geq n\}$ is a
stopping time with respect to the natural filtration of
$(Z_{\mathbf{i}})_{\vert\mathbf{i} \vert\leq n}$,
strong Markov property gives
\[
\mathbb{P}(L=0)\leq\mathbb{P}(T_n<\infty) \mathbb{P}(L=0)^n +
\mathbb{P}(T_n=\infty).
\]
If $\mathbb{P}(L=0)<1$, letting $n\rightarrow\infty$ gives
\[
\mathbb{P}(L=0)\leq\lim_{n\rightarrow\infty} \mathbb{P}(T_n=\infty)
=\mathbb{P}(\mathcal{Z}_n \mbox{ is bounded})=\mathbb{P}(\operatorname{Ext})
\]
since $\mathcal{Z}_n$ is a BGW. This completes the proof. One can also
use a coupling argument: the number of contaminated cells
starting with one single cell with $n$ parasites is less than the
number of contaminated cells
starting from $n$ cells with one single parasite.
\end{pf*}

\section{Tree of contaminated cells}\label{s4}

Here, we prove that contaminated cells are not concentrated in a cell
line. If
$m_0+m_1\leq1$, contaminated
cells die out but conditionally on the survival of parasites at
generation $n$, the number of leaves of the tree of contaminated
cells tends to $\infty$ as $n \rightarrow\infty$. The proof of this
result will also ensure that,
if $m_0+m_1>1$, the number of contaminated cells tends to $\infty$
provided that they do not die out.
\begin{Thm}\label{nb}
If $m_0+m_1\leq1$, $\# \{\mathbf{i} \in\mathbb{T}\dvtx Z_{\mathbf{i}}\ne0,
Z_{\mathbf{i}0}=0, Z_{\mathbf{i}1}=0\}$ conditioned by $\# \mathbb{G}^*_n >0$
converges in probability as $n\rightarrow\infty$ to $\infty$.

If $m_0+m_1>1$, conditionally on $\operatorname{Ext}^c$, $\# \mathbb{G}
^*_n\stackrel{n\rightarrow\infty}{\longrightarrow}\infty$ a.s.
\end{Thm}
\begin{Rque*}
In the conditions of the theorem, $\# \mathbb{G}^*_n$ (resp. the number
of leaves) grows at least linearly with respect to $n$ (see Section
\ref{s5} for further results).
In the case $m_0+m_1\leq1$, conditionally on $\# \mathbb{G}^*_n >0$,
the tree of contaminated cells is a spine with finite
subtrees, as for BGW conditioned to survive (see \cite{Geig,Lyons}).
\end{Rque*}

We need two lemmas for the proof. First we prove that the ancestor of a
contaminated cell has given birth to two contaminated cells
with a probability bounded from below. We have to distinguish the case
where $\mathbb{P}(Z^{(0)}Z^{(1)}=0)=1$, since in that case a cell must
contain at
least two parasites so that it can give birth to two contaminated cells.
\begin{Lem} \label{nbone}
There exists $\alpha>0$ such that for all $N \in\mathbb{N}$,
$\mathbf{i}\in\mathbb{G}_N$, $n<N$ and $k\geq2$,
\[
\mathbb{P}( Z_{\mathbf{j}0}\ne0, Z_{\mathbf{j}1}\ne0  \mid
Z_{\mathbf{j}}=k, Z_{\mathbf{i}}>0)\geq\alpha
\]
denoting $\mathbf{j}=\mathbf{i}\mid n$. If $\mathbb{P}
(Z^{(0)}Z^{(1)}=0)\ne1$, this result also holds for $k=1$.
\end{Lem}
\begin{pf}
We consider first the case $\mathbb{P}(Z^{(0)}Z^{(1)}=0)\ne1$ and we choose
$(k_{0},k_{1}) \in\mathbb{N}^{*2}$ such that
$\mathbb{P}( (Z^{(0)},Z^{(1)})=(k_0,k_{1}))>0$. For every $k\in
\mathbb{N}^*$, we have
\[
\mathbb{P}( Z_{\mathbf{j}0}\ne0, Z_{\mathbf{j}1}\ne0  \mid
Z_{\mathbf{j}}=k, Z_{\mathbf{i}}>0 )
\geq \mathbb{P}( Z_{\mathbf{j}0}\ne0, Z_{\mathbf{j}1}\ne0  \mid
Z_{\mathbf{j}}=1, Z_{\mathbf{i}}>0 ).
\]
%
Moreover, as the function $\mathbb{R}^*_+\ni u \mapsto
\frac{1-e^{-u}}{u}$ decreases, we have for all $y,x> 0$ and $p \in
[0,1[$,
%
\begin{equation}\label{minunif}
\frac{1-p^{x}}{1-p^{y}}\geq\frac{x}{\max \{y,x\}}.
\end{equation}
Let $a\in\{0,1\}$ and $\mathbf{k}$ such
that $\mathbf{i}=\mathbf{j}a\mathbf{k}$. Then for all
$(k'_{0},k'_{1}) \in
\mathbb{N}^{2}-(0,0)$,
\begin{eqnarray*}
&& \frac{ \mathbb{P}(Z_{\mathbf{j}0}= k_{0}, Z_{\mathbf{j}1}=k_{1}
\mid  Z_{\mathbf{j}}=1, Z_{\mathbf{i}}>0 )}{\mathbb{P}(Z_{\mathbf{j}0}= k'_{0},
Z_{\mathbf{j}1}= k'_{1}  \mid   Z_{\mathbf{j}}=1, Z_{\mathbf{i}}>0 )}
\\
&& \qquad=  \frac{ \mathbb{P}(Z^{(0)}= k_{0}, Z^{(1)}= k_{1}  \mid
Z_{a\mathbf{k}}>0 )}{\mathbb{P}(Z^{(0)}= k'_{0}, Z^{(1)}= k'_{1}  \mid   Z_{a\mathbf{k}}>0 )}
\\
&& \qquad =  \frac{\mathbb{P}(Z_{a\mathbf{k}}>0  \mid Z^{(0)}=k_0,
Z^{(1)}=k_1) }{ \mathbb{P}(Z_{a\mathbf{k}}>0  \mid Z^{(0)}=k'_0, Z^{(1)}=k'_1)}
\frac{\mathbb{P}( Z^{(0)}= k_{0}, Z^{(1)}= k_{1} )}{\mathbb{P}(Z^{(0)}= k'_{0},
Z^{(1)}= k'_{1} )}
\\
&& \qquad=   \frac{1-\mathbb{P}(Z_{\mathbf{k}}=0)^{k_a}}{1-\mathbb{P}
(Z_{\mathbf{k}}=0)^{k'_a}}\frac{\mathbb{P}( (Z^{(0)},Z^{(1)})=(k_0,k_{1}))}
{\mathbb{P}((Z^{(0)},Z^{(1)})=(k'_{0},k'_{1}))}
\\
&& \qquad\geq  \frac{\min  \{k_0,k_1 \}}{k_0+k_1+k_0'+k'_1}
\frac{\mathbb{P}( (Z^{(0)},Z^{(1)})=(k_0,k_{1}))}
{\mathbb{P}((Z^{(0)},Z^{(1)})=(k'_{0},k'_{1}))} \qquad\mbox{using
(\ref{minunif})}.
\end{eqnarray*}
Cross product and sum over $(k'_{0},k'_{1})$ give
\begin{eqnarray*}
&& \bigl [\mathbb{E}\bigl(Z^{(0)}+Z^{(1)}\bigr)+k_0+k_1 \bigr]
\mathbb{P}(Z_{\mathbf{j}0}= k_{0},
Z_{\mathbf{j}1}= k_{1} \mid Z_{\mathbf{j}}=1, Z_{\mathbf{i}}>0)
\\
&&\qquad \geq\min \{k_0,k_1 \} \mathbb{P}
\bigl( \bigl(Z^{(0)},Z^{(1)}\bigr)=(k_0,k_{1})\bigr).
\end{eqnarray*}
This gives the result since $\mathbb{P}(Z_{\mathbf{j}0}= k_{0},
Z_{\mathbf{j}1}= k_{1} \mid Z_{\mathbf{j}}=1, Z_{\mathbf{i}}>0)\geq\alpha$ with
\[
\alpha=\frac{\min  \{k_0,k_1 \} \mathbb{P}( (Z^{(0)},
Z^{(1)})=(k_0,k_{1})) }{\mathbb{E}(Z^{(0)}+Z^{(1)})+k_0+k_1}>0.
\]

In the case $\mathbb{P}(Z^{(0)}Z^{(1)}=0)=1$, we choose $(k_{0},k_{1})
\in\mathbb{N}^{*2}$ such that $\mathbb{P}_2( (Z_{0},\break Z_{1})=(k_0,k_{1}))>0$ [using
(\ref{condp})]. We make then the same
proof as above with $Z_{\mathbf{j}}=2$ and
\[
\alpha=\frac{\min \{k_0,k_1 \} \mathbb{P}_2( (Z_{0},Z_1)=( k_{0},
k_{1})) }{\mathbb{E}_2(Z_{0}+Z_{1})+k_0+k_1},
\]
so that the result follows as previously.
\end{pf}

Thus if $\mathbb{P}(Z^{(0)}Z^{(1)}=0)=1$, we need to prove that there
are many
cells with more than two parasites in a
contaminated cell line.
\begin{Lem}\label{nbtwo}
If $\beta:=\mathbb{P}( Z^{(0)} \geq2 \mbox{ or } Z^{(1)}\geq2)>0$ then
\[
\inf_{\mathbf{i} \in\mathbb{G}_n} \mathbb{P}(\# \{ \mathbf
{j}<\mathbf{i}\dvtx Z_{\mathbf{j}0}\geq2
\mathrm {or} \ Z_{\mathbf{j}1}\geq2 \}\geq\beta n/2  \vert
Z_{\mathbf{i}}>0)\stackrel{n\rightarrow\infty}{\longrightarrow}1.
\]
\end{Lem}
\begin{pf}
For all $\mathbf{i}\in\mathbb{G}_n$ and $\mathbf{j}<\mathbf{i}$,
let $\mathbf{k}$ such that $\mathbf{i}=\mathbf{j}\mathbf{k}$,
then for every $\alpha>0$,
\[
\mathbb{P}( Z_{\mathbf{j}0} \geq2 \mbox{ or }
Z_{\mathbf{j}1}\geq2  \mid  Z_{\mathbf{j}}=\alpha,
Z_{\mathbf{i}}>0 )\geq \mathbb{P}( Z_{0} \geq2 \mbox{ or }
Z_{1}\geq2  \mid Z_{\mathbf{k}}>0) \geq \beta.
\]
Then conditionally on $ Z_{\mathbf{i}}>0$, $\# \{
\mathbf{j}<\mathbf{i} \dvtx Z_{\mathbf{j}0}\geq2$ or
$Z_{\mathbf{j}1}\geq2 \}\geq \sum_{k=0}^n \beta_k$, where
$(\beta_k)_{1\leq k\leq n}$ are i.i.d.
and distributed as a Bernoulli($\beta$). Conclude with the law of
large numbers.
\end{pf}
\begin{pf*}{Proof of Theorem \ref{nb}}
We consider first the case when $m_0+m_1>1$, work conditionally on
$\operatorname{Ext}^c$ and choose $\mathbf{i}\in\delta\mathbb{T}^*$.

If $\mathbb{P}(Z^{(0)}Z^{(1)}=0)\ne1$, Lemma \ref{nbone} (with $k=1$)
entails that a.s. under \mbox{$\mathbb{P}(\cdot\vert Z_{\mathbf{i}}>0)$},
\[
\# \{\mathbf{j}<\mathbf{i}\dvtx Z_{\mathbf{j}0}>0, Z_{\mathbf{j}1}>0 \}=\infty.
\]
Using the branching property and the fact that the probability of
nonextinction of parasites is positive ensures that
$\# \mathbb{G}^*_n \stackrel{n\rightarrow\infty}{\longrightarrow}\infty
$ a.s.

If $\mathbb{P}(Z^{(0)}Z^{(1)}=0)=1$ then $\mathbb{P}( Z^{(0)} \geq2$
or $Z^{(1)}\geq2 )>0$ and by Lemma \ref{nbtwo}, we have a.s. on
$\mathbb{P}(\cdot\vert Z_{\mathbf{i}}>0)$,
\[
\# \{\mathbf{j}<\mathbf{i} \dvtx Z_{\mathbf{j}0}\geq2
\mbox{ or } Z_{\mathbf{j}1}\geq2 \}=\infty.
\]
Using as above Lemma \ref{nbone} (with $k=2$) and the fact that the
probability of nonextinction of parasites is positive ensures that
$\# \mathbb{G}^*_n \stackrel{n\rightarrow\infty}{\longrightarrow}\infty$ a.s.

We consider now the case when $m_0+m_1\leq1$ and work conditionally on
$\mathbf{i}=(\alpha_1,\ldots,\alpha_n)\in\mathbb{G}_n^*$. We denote
$\mathbf{i}_j:=(\alpha_1,\ldots,\alpha_{j-1},1-\alpha_j)$ for
$1\leq j \leq n$.

If $\mathbb{P}(Z^{(0)}Z^{(1)}=0)\ne1$, Lemma \ref{nbone} entails that
%
\begin{equation}\label{prcond}
\forall1\leq j \leq n, k\geq1 \qquad
\mathbb{P}(Z_{\mathbf{i}_j}>0 \mid Z_{\mathbf{i}\vert j-1}=k,
Z_{\mathbf{i}}>0)\geq \alpha.
\end{equation}
Moreover, if $Z_{\mathbf{i}_j}>0$, then the tree of
contaminated cells
rooted in $\mathbf{i}_j$
dies out and so has at least one leaf. So by the branching property,
the number of
leaves converges in probability to infinity as $n$ tends to
infinity.

If $\mathbb{P}(Z^{(0)}Z^{(1)}=0)=1$, (\ref{prcond}) holds for $k\geq2$ and
Lemma \ref{nbtwo} allows to conclude similarly in this case.
\end{pf*}

\section{Proportion of contaminated cells with a given number of
parasites}\label{s5}

We determine here the asymptotics of the number of contaminated cells
and the proportion
$F_k$ of cells with $k$ parasites, defined as
\[
F_k(n):= \frac{ \# \{ \mathbf{i}\in\mathbb{G}_n ^*\dvtx
Z_{\mathbf{i}}=k\}}{\# \mathbb{G}_n^*} \qquad (k\in\mathbb{N}^*).
\]
In that view,
we introduce the Banach space $l^1(\mathbb{N})$ and the subset of
frequencies $\mathbb{S}^1(\mathbb{N})$ which we endow with the norm
$\|\cdot \|_1$ defined by:
\begin{eqnarray*}
l^1(\mathbb{N}) &:= & \Biggl\{(x_i)_{i\in\mathbb{N}}\dvtx
\sum_{i=0}^{\infty}\vert x_i\vert <\infty\Biggr\}, \qquad
\|(x_i)_{i\in\mathbb{N}} \|_1=\sum_{i=0}^{\infty}\vert x_i\vert,
\\
\mathbb{S}^1(\mathbb{N}) &:=&
\Biggl\{(f_i)_{i\in\mathbb{N}}\dvtx
\forall i \in \mathbb{N},  f_i \in\mathbb{R} ^+,
\sum_{i=0}^{\infty} f_i=1\Biggr\}.
\end{eqnarray*}
We shall work conditionally on $\operatorname{Ext}^c$ or $\mathcal{Z}_n>0$
and introduce
%
\begin{equation}\label{Pstar}
\mathbb{P}^*:=\mathbb{P}( \cdot \mid\operatorname{Ext}^c), \qquad
\mathbb{P}^n:=\mathbb{P}( \cdot  \mid \mathcal{Z}_n>0).
\end{equation}

The asymptotics of the proportions depend naturally on the
distribution of ($Z^{(0)},Z^{(1)}$) and we determine five
different behaviors according to the bivariate value of $(m_0,m_1)$.

The proofs of the convergences use the asymptotic distribution
of the number of parasites of a typical contaminated cell at
generation $n$, which is equal to \mbox{$\mathbb{P}^n(Z_{\mathbf{U}_n}\in
\cdot)$}, where
$\mathbf{U}_n$ is a uniform random variable in $\mathbb{G}_n^*$ independent
of $(Z_{\mathbf{i}})_{\mathbf{i}\in\mathbb{T}^*}$. This
distribution is different
from the distribution of $Z_n^*$, that is the number of parasites
of a random cell line conditioned to be contaminated at
generation~$n$. The following example even proves that
$\mathbb{P}^n(Z_{\mathbf{U}_n}\in\cdot)$ and $\mathbb{P}(Z_n^*\in\cdot)$
could be a priori very different.
\begin{ex*}\label{choixunif}
Suppose that generation $n$ (fixed) contains
$100$ cells with $1$ parasite (and no other contaminated cells) with
probability $1/2$ and it contains $1$
cell with $100$ parasites with probability $1/2$ (and no other
contaminated cells). Compare then
\begin{eqnarray*}
\mathbb{P}^n(Z_{\mathbf{U}_n}=1) &=& 1/2, \qquad\hspace*{11pt}
\mathbb{P}^n(Z_{\mathbf{U}_n}=100)=1/2;
\\
\mathbb{P}(Z^*_{n} = 1) &=& 100/101, \qquad
\mathbb{P}(Z^*_{n} = 100)=1/101.
\end{eqnarray*}
\end{ex*}

Actually the convergence of $(Z_n^*)_{n\in\mathbb{N}}$ leads to the result
obtained by Kimmel \cite{kim}
in the continuous analogue of this model. That is,
%
\begin{equation}\label{kim}
\frac{\mathbb{P}(Z_n=k)}{\mathbb{P}(Z_n>0)}= \frac{\sum_{\mathbf
{i}\in\mathbb{G}_n}
\mathbb{E}(\mathbh{1}_{Z_{\mathbf{i}}=k})}{\sum_{\mathbf{i}\in
\mathbb{G}_n} \mathbb{E}
(\mathbh{1}_{Z_{\mathbf{i}}>0})}=\frac{\mathbb{E}(\# \{ \mathbf
{i}\in\mathbb{G}_n \dvtx
Z_{\mathbf{i}}=k\})}{\mathbb{E}(\# \mathbb{G}_n^*)}
\end{equation}
tends to $\mathbb{P}(\Upsilon=k)$ whereas we are here interested in the
expectation of $F_k(n)$.

A sufficient condition to get the equality of the two distributions is
that $\# \mathbb{G}_n^*$
is deterministic, which does not hold here.
But in the case when $(m_0,m_1)\in D_3$, we shall prove that
$\#\mathbb{G}_n^*$ is asymptotically
proportional to $(m_0+m_1)^n$ as \mbox{$n \rightarrow\infty$} (forthcoming
Proposition \ref{estim}). This
enables us to control $\mathbb{P}^n(Z_{\mathbf{U}_n}\in\cdot)$ by the
distribution of $\mathbb{P}(Z_n^*\in\cdot )$. More precisely,
it is sufficient to prove the separation of descendances of parasites
(Proposition \ref{separation}) and the control of filled-in cells
(Lemma~\ref{contr}) using the results about the BPRE $Z_n^*$.
These two results are the keys for Theorems
\ref{limf1}, \ref{limf4} and \ref{limf2}. Similarly, when
$(m_0,m_1)\in D_5$, we already know that $\# \mathbb{G}_n^*$ is
approximatively
equivalent to $2^n$. Then the fact that
$Z_n^*$ explodes as $n\rightarrow\infty$ (by Proposition
\ref{explo}) will ensure that the proportion of filled-in cells
among contaminated cells tends to one
(Theorem \ref{limf3} below).

\subsection{Case $(m_0,m_1) \in D_5$ \textup{(}$m>1$\textup{)}}\label{s5.1}

In that case, recall that conditionally
on $\operatorname{Ext}^c$, $\#\mathbb{G}_n^*$ is asymptotically
proportional to $2^n$
(by Theorem \ref{Text}).
Moreover the contaminated cells become largely infected, as stated below.
\begin{Thm}\label{limf3}
Conditionally on $\operatorname{Ext}^c$, for every
$k\in\mathbb{N}$, $F_k(n)$ converges in probability to $0$ as
$n\rightarrow \infty$, that is,
\[
\forall K,\varepsilon>0 \qquad\mathbb{P}^*
\biggl( \frac{\# \{ \mathbf{i}\in\mathbb{G}_n \dvtx
Z_{\mathbf{i}}\geq K \}}{\# \mathbb{G}_n^*}\geq1-\varepsilon \biggr)
\stackrel{n\rightarrow\infty}{\longrightarrow} 1.
\]
\end{Thm}

If $m_0=m_1$, the number of parasites in a contaminated cell is of the
same order as $m_0^n$. More precisely, for every $\varepsilon>0$,
\[
\sup_{n\in\mathbb{N}} \biggl\{\mathbb{P}^* \biggl(\frac{\#\{ \mathbf{i} \in\mathbb{G}^*_n \dvtx
Z_{\mathbf{i}}\leq\alpha m_0^n\}}{\# \mathbb{G}_n^*}\geq \varepsilon \biggr) \biggr\}
\stackrel{\alpha\rightarrow0}{\longrightarrow}0.
\]
\begin{pf}
In that case, use Theorem \ref{Text} and (\ref{Pstar}) to get
that there exists a nonnegative random variable $\widetilde {L}$ such
that
%
\begin{equation}\label{Lpos}
\# \mathbb{G}_n^*\geq2^n \widetilde{L}, \qquad
\mathbb{P}^*(\widetilde {L}=0)=0.
\end{equation}
Let $K,\eta$ and $\varepsilon>0$ and put
$B_n(K,\eta):= \{\frac{\# \{ \mathbf{i} \in\mathbb{G}_n^* \dvtx
Z_{\mathbf{i}}\leq K\} }{\# \mathbb{G}^*_n}\geq\eta \}\cap\operatorname{Ext}^c$, then
\[
\sum_{\mathbf{i} \in\mathbb{G}_n^*} \mathbh{1}_{\{Z_{\mathbf{i}}\leq K\}}
\geq\eta2^n \widetilde {L} \mathbh{1}_{B_n(K,\eta)}
\]
which gives, taking expectations,
\[
\mathbb{E}\bigl(\widetilde {L}\mathbh{1}_{B_n(K,\eta)}\bigr)
\leq\frac{\mathbb{E}(\sum_{\mathbf{i} \in\mathbb{G}_n^*} 2^{-n}\mathbh{1}_{\{Z_{\mathbf{i}}\leq K\}})}
{\eta}=\frac{\mathbb{P}(0<Z_n\leq K)}{\eta}.
\]
Use then Proposition \ref{explo} and (\ref{Lpos}) to choose $n$ large
enough so that
\[
\mathbb{P}(B_n(K,\eta))\leq\varepsilon,
\]
which completes the proof of the theorem. In the case $m_0=m_1=m$,
follow the proof above
and use that $Z_n/m^n$ converges to a positive limit on $\operatorname{Ext}^c$
(see \cite{atr}) to get the finer result given after the theorem.
\end{pf}
%

\subsection{Case $(m_0,m_1) \in D_3$ \textup{(}$m\leq1$\textup{)}}\label{s5.2}

We assume here $\mathbb{E}(Z^{(a)2})<\infty$ and prove that $(F_k(n))_{k
\in\mathbb{N}}$ converges to a deterministic limit. We prove the
convergence thanks to the Cauchy criterion [using completeness of
$l^1(\mathbb{N})$]. The fact that the limit is deterministic is a
consequence of the separation of the descendances of parasites and
the law of large numbers. Once we know this limit is
deterministic, we identify it with the Yaglom limit $\Upsilon $ (see
Section~\ref{s6.1} for proofs).
\begin{Thm}\label{limf1}
Conditionally on $\operatorname{Ext}^c$, as $n \rightarrow
\infty$, $(F_k(n))_{k \in\mathbb{N}}$ converges in probability in
$\mathbb{S}^1(\mathbb{N})$ to ($\mathbb{P}(\Upsilon =k))_{k \in
\mathbb{N}}$.
\end{Thm}
\begin{Rque*} \label{rem}
We get here a realization of the Yaglom
distribution $\Upsilon $.

The limit just depends on the one-dimensional distributions of
$(Z^{(0)},Z^{(1)})$. More precisely,
recall that the probability generating function $G$ of $\Upsilon $ is
characterized by (\ref{equatfonct}).

This theorem still holds starting from $k$ parasites. We also
easily get a similar result
in the case when a cell gives birth to $N$ cells ($N\in\mathbb{N}$).

As an application, we can obtain numerically
the Yaglom quasistationary distribution of any BGW. Let $Z$ be
the reproduction law of a BGW with mean $m<1$ and choose $N$ such
that $Nm>1$. Consider Kimmel's model where each cell divides into
$N$ daughter cells and
$Z^{(0)}\stackrel{d}{=}Z^{(1)}\stackrel{d}{=}\cdots\stackrel{d}{=}
Z^{(N)}\stackrel{d}{=}Z$. Computing then the asymptotic of the
proportions of contaminated cells with $k$ parasites gives the
Yaglom quasistationary distribution associated to $Z$. If
$\mathbb{P}(\operatorname{Ext})\ne0$, one can
start from many parasites ``to avoid'' extinction.

More generally, we can obtain similarly the Yaglom quasistationary distribution
of any BPRE with finite number $k$ of environments such that
$\sum_{1}^{k} m_i^2<\sum_{1}^{k} m_i$.

This theorem is in the same vein as Theorem 11 in \cite{guy}. But we
can not follow the same approach
as Guyon for the proof. Indeed we have to consider here the
proportions among the contaminated cells in generation $n$ whereas
Guyon considers proportions among all cells in generation $n$.
Unfortunately, the subtree of contaminated cells is itself random
and induces long-range dependences between cells lines, so that
Guyon's arguments do not hold here. Moreover, Theorem 11 in
\cite{guy} relies on an ergodicity hypothesis which cannot be
circumvented.
\end{Rque*}
\begin{ex*}
We give two examples when the limit can be calculated.

Trivial case: $\mathbb{P}(Z^{(0)} \in\{0,1\}, Z^{(1)} \in\{0,1\})=1$ leads to
$\mathbb{P}(\Upsilon=1)=1$.

Symmetric linear fractional case: $p \in\,]0,1[$,
$b \in\,]0,(1-p)^2[$ and
\[
\mathbb{P}\bigl(Z^{(0)}=k\bigr)
= \mathbb{P}\bigl(Z^{(1)}=k\bigr)=bp^{k-1} \qquad \mbox{if } k\geq1
\]
and $\mathbb{P}(Z^{(0)}=0)=\mathbb{P}(Z^{(1)}=0)=(1-b-p)/(1-p)$.
Then $m_0=m_1=b/(1-p)^2<1$ and letting
$s_0$ be the root of $f_0(s)=s$ larger than $1$,
\[
\forall k\geq1 \qquad\mathbb{P}(\Upsilon=k)=(s_0-1)/s_0^k.
\]
\end{ex*}

As asymptotically we know the number of parasites and the
proportion of cells with $k$ parasites, we get the number of
contaminated cells [recall that $W$ is given by (\ref{nbparas})].
\begin{Cor}\label{gn} 
Conditionally on $\operatorname{Ext}^c$, the following
convergences hold in probability
\[
\frac{\# \mathbb{G}^*_n}{\mathcal{Z}_n} \stackrel{n\rightarrow\infty}
{\longrightarrow}\frac{1}{\mathbb{E}(\Upsilon )},
\qquad
\frac{\# \mathbb{G}^*_n}{(m_0+m_1)^n}
\stackrel{n\rightarrow\infty}{\longrightarrow} \frac{W}{\mathbb{E}(\Upsilon )}.
\]
\end{Cor}

We can also consider the ancestors at generation $n$ of the cells of
$\partial\mathbb{T}^*$, which amounts to considering
\[
F_k(n,p)=\frac{\# \{ \mathbf{i} \in\mathbb{G}_{n+p}^* \dvtx
Z_{\mathbf{i}\mid n}=k\}}{\# \mathbb{G}_{n+p}^* }
\]
and let $p\rightarrow\infty$. Letting then $n \rightarrow\infty$
yields the biased Yaglom
quasistationary distribution, thanks to the separation of
descendances of parasites.
\begin{Cor}\label{limf13}
Conditionally on $\operatorname{Ext}^c$, for every
$k\in\mathbb{N}$, $F_k(n,p)$ converges in probability in
$\mathbb{S}^1(\mathbb{N})$ as
$p$ tends to infinity. This limit converges in probability in
$\mathbb{S}^1(\mathbb{N})$ as $n \rightarrow\infty$:
\[
\forall k \in\mathbb{N} \qquad
\lim_{n\rightarrow\infty} \lim_{p\rightarrow\infty} F_k(n,p)
\stackrel{\mathbb{P}}{=} \frac{k\mathbb{P}(\Upsilon =k)}{\mathbb{E}(\Upsilon )}.
\]
\end{Cor}

We get here an interpretation of the fact that the stationary
distribution of the Q-process associated to the BPRE
$(Z_n)_{n\in\mathbb{N}}$ is the size-biased Yaglom limit (see~\cite{afa2}).

\subsection{Case $(m_0,m_1) \in D_2$}\label{s5.3}

In that case, the parasites die out. So we condition by $\mathcal{Z}_n>0$,
we still assume $\mathbb{E}(Z^{(a) \ 2})<\infty$ and we get a similar
result.
\begin{Thm}\label{limf4}
As $n \rightarrow\infty$, $(F_k(n))_{k \in\mathbb{N}}$ conditioned
by $\mathcal{Z}_n>0$ converges in distribution
on $\mathbb{S}^1(\mathbb{N})$ to ($\mathbb{P}(\Upsilon =k))_{k \in\mathbb{N}}$.
\end{Thm}

The proof follows that of the previous theorem. Indeed
$(\ref{conH})$ is still satisfied and we can use the same results
on the BPRE $(Z_n)_{n\in\mathbb{N}}$. There are only two differences.
First, we work under $\mathbb{P}^n$ instead of $\mathbb{P}^*$.
Moreover $\mathcal{Z}_n$
satisfies now $\mathbb{P}(\mathcal{Z}_n>0)\stackrel{n\rightarrow\infty}{\sim}
2/(\mathrm {Var}(Z^{(0)}+Z^{(1)})n)$ and $\mathcal{Z}_n/n$
conditioned to be nonzero converges in distribution as $n \rightarrow\infty$ to an
exponential variable $\mathcal{E}$ of parameter $2/(\hat{m}+1)$
(see Section~\ref{s2.1}). As above, we can derive the following result.
\begin{Cor} \label{ggn}
As $n \rightarrow\infty$, $\# \mathbb{G}_n^*/n$ conditioned by
$\# \mathbb{G}_n^*>0$ converges in distribution to
$\mathcal{E}/\mathbb{E}(\Upsilon )$.
\end{Cor}
%

\subsection{Case $(m_0,m_1) \in D_1$}\label{s5.4}

In this case, the number of contaminated cells does not explode
and the number of cells of type $k$ at generation $n$ conditioned
by the survival of parasites in this generation converges weakly
to a nondeterministic limit (see Section \ref{s7} for proofs).
\begin{Thm}\label{limf2}
As $n \rightarrow\infty$, $(\#\{\mathbf{i} \in\mathbb{G}_n^* \dvtx
Z_{\mathbf{i}}=k \})_{k \in\mathbb{N}}$
conditioned on $\mathcal{Z}_n>0$ converges in distribution on
$l^1(\mathbb{N})$
to a random sequence ($N_k)_{k \in\mathbb{N}}$ which satisfies
$\mathbb{E}(\sum_{k\in\mathbb{N}} kN_{k})<\infty$.
\end{Thm}

As above, we get:
\begin{Cor} \label{gggn}
$\# \mathbb{G}^*_n$ conditioned by $\# \mathbb{G}^*_n >0$ converges
in distribution to a positive finite random variable.
\end{Cor}

Picking a cell uniformly on $\partial\mathbb{T}^*$ leads again to the
size-biased distribution.
\begin{Cor}\label{limf22}
For every $n\in\mathbb{N}$, $(\# \{ \mathbf{i} \in\mathbb{G}
_{n+p}^* \dvtx Z_{\mathbf{i}\mid n}=k\})_{k \in\mathbb{N}}$
conditioned on $\mathcal{Z}_{n+p}>0$ converges weakly in
$l^1(\mathbb{N})$ to a random sequence as $p$ tends to infinity. This
limit converges weakly as $n \rightarrow\infty$.
\[
\forall k \in\mathbb{N} \qquad
\lim_{n\rightarrow\infty}\lim_{p\rightarrow\infty}\#
\{ \mathbf{i} \in\mathbb{G}_{n+p}^* \dvtx
Z_{\mathbf{i}\mid n}=k\} \vert\mathcal{Z}_{n+p}>0
= \frac{kN_k}{\sum_{k'\in\mathbb{N}} k'N_{k'}}.
\]
\end{Cor}
%

\subsection{Remaining domain\textup{:} $(m_0,m_1) \in D_4$}\label{s5.5}

In this domain, the asymptotic of the mean of the number of
contaminated cells, that is $\mathbb{E}(\# \mathbb{G}_n^*)=2^n\mathbb{P}(Z_n>0)$, is
different from the previous ones.

Recalling Section \ref{s2.2}, this asymptotic depends on three
subdomains, the interior of $D_4$ and the two connex components of
its boundary. More precisely, it depends on $m_0m_1=1$ or
$m_0m_1<1$ and $m_0\log (m_0)+m_1\log (m_1)$ is
positive or
zero.

If $(m_0,m_1)\in D_4$ and $m_0<1<m_1$, using (\ref{esp}) and a
coupling argument with Corollary \ref{gn}, one can prove that
\[
\sup_{n\in\mathbb{N}}
\biggl\{\mathbb{P} \biggl(\frac{\# \mathbb{G}_n^*}{2^n\mathbb{P}(Z_n>0)}
\geq A,  \frac{\# \mathbb{G}_n^*}{(m_0+\widetilde {m_0})^n}\leq1/A \biggr)
\biggr\}\stackrel{A\rightarrow 0}{\longrightarrow}0,
\]
where
$\widetilde {m_0}= (1+\sqrt{1+4(m_0-m_0^2)} )/2>1$. Thus
$\#\mathbb{G}_n^*$
grows geometrically and one can naturally conjecture that
$\#\mathbb{G}_n^*$ is asymptotically proportional to
$\mathbb{E}(\#\mathbb{G}_n^*)=2^n\mathbb{P}(Z_n>0)$.

Moreover separation of descendances of parasites, control of
filled-in cells and Corollary \ref{limf13} do not hold in this
case. Thus determining the limit behaviors here requires a
different approach.

Finally, note that in the subdomain $m_0m_1=1$ (boundary of
$D_5$), $(Z_n^*)_{n\in\mathbb{N}}$ explodes (see \cite{afa}) so the
asymptotic proportion of contaminated cells which are arbitrarily
largely contaminated should be equal to $1$ as in Theorem
\ref{limf3}.

%
%

\section{Proofs in the case $(m_0,m_1)\in D_3$}\label{s6}

We assume in this section that $\mathbb{E}(Z^{(a)2})<\infty$ (i.e.,
$\widetilde {m}<\infty$) and we start with giving
some technical results.

\subsection{Preliminaries}\label{s6.1}

First, note that for all $u,v$ $\in$ $l^1(\mathbb{N}^*)$, we have
%
\begin{equation}\label{impliquer}
\biggl\| \frac{u}{\| u \|_1}-\frac{v}{\| v\| _1}\biggr\|_1
= \biggl\| \frac{u-v}{\| u\| _1}+ \frac{v}{\| v\| _1}
\frac{\| v\|_1-\| u\| _1}{\| u\| _1}\biggr\|_1
\leq2\frac{\| u-v\| _1}{\| u\| _1}.
\end{equation}

Moreover by (\ref{nbparas}), there exist two random variables $C$
and $D$ a.s finite such that
%
\begin{equation}\label{minW}
\qquad
\forall n \in\mathbb{N} \qquad
C\leq\frac{\mathcal{Z}_n}{(2m)^n}\leq D \qquad\mbox{a.s.},\
\mathbb{P}^*(C=0)=\mathbb{P}^*(D=0)=0
\end{equation}
and as $\bigcap_{n\in\mathbb{N}}
\{\mathcal{Z}_n>0\}=\{\forall n \in\mathbb{N}\dvtx \mathcal{Z}_n>0\}$,
we have
%
\begin{equation}\label{P_n}
\sup_{A}\{ \vert\mathbb{P}^n(A)-\mathbb{P}^*(A) \vert\}
\stackrel{n\rightarrow\infty}{\longrightarrow}0.
\end{equation}

We focus now on the BPRE $(Z_n)_{n\in\mathbb{N}}$. First, by
induction and
convexity of $f_a$, we have for every $\mathbf{i} \in\mathbb{G}_n$ (see
Section \ref{s2.2} for the notation)
%
\begin{equation}\label{majsur}
\mathbb{P}(Z_{\mathbf{i}}>0)=1-f_{\mathbf{i}}(0)\leq m_{\mathbf{i}}.
\end{equation}
Then identities
(\ref{majsur}) and (\ref{equiv}) entail that there exists $M>0$
such that
%
\begin{equation}\label{majP}
M\leq\frac{\mathbb{P}(Z_n>0)}{m^n}\leq 1.
\end{equation}

Moreover, by Corolary 2.3 in \cite{afa2}, we have
%
\begin{equation}\label{control}
\lim_{K\rightarrow\infty} \sup_{n\in\mathbb{N}}
\{ \mathbb{E}( Z_n \mathbh{1}_{Z_n \geq K} \mid Z_n>0) \} =0.
\end{equation}

Finally, following the proof of Theorem 1.2 in \cite{Guiv} (see
\cite{Sub}, Section 2.1 for details) ensures that, if
$(Z_n^{(1)})_{n\in\mathbb{N}}$ and $(Z_n^{(2)})_{n\in\mathbb{N}}$
are two
independent BPRE distributed as $(Z_n)_{n\in\mathbb{N}}$, we have
\[
\mathbb{P}\bigl(Z_n^{(1)}>0, Z_n^{(2)}>0\bigr)
= o\bigl(\mathbb{P}(Z_n>0)\bigr)=o(m^n) \qquad
(n\rightarrow\infty).
\]
Then, we have
%
\begin{equation}\label{negli}
2^{-n} \sum_{\mathbf{i}\in\mathbb{G}_n}\mathbb{P}(Z_{\mathbf{i}}>0)^2=o(m^n)
\qquad(n\rightarrow\infty).
\end{equation}

\subsection{Estimation of $ \# \mathbb{G}_n^*$}\label{s6.2}

We prove here that the number of parasites which belong to filled-in
cells is
negligible compared to the total number of parasites (see also
Lemma \ref{contr} for a result of the same kind). To prove this result,
we use its counterpart for BPRE $(Z_n)_{n\in\mathbb{N}}$ conditioned
to be nonzero.
\begin{Lem}\label{celltc}
For every $\eta>0$,
\[
\sup_{n \in\mathbb{N}} \biggl\{\mathbb{P}^*
\biggl( \frac{ \sum_{\mathbf{i} \in\mathbb{G}_n^*} Z_{\mathbf{i}}
\mathbh{1}_{\{Z_{\mathbf{i}}> K\}}}{\mathcal{Z}_n} \geq\eta
\biggr) \biggr\}\stackrel{K\rightarrow\infty}{\longrightarrow}0.
\]
\end{Lem}
\begin{pf}
Let $\eta>0$ and write
\[
A_n(K,\eta):= \biggl\{ \frac{ \sum_{\mathbf{i}
\in\mathbb{G}_n^*} Z_{\mathbf{i}} \mathbh{1}_{\{Z_{\mathbf{i}}> K\}}}{\mathcal{Z}_n}
\geq\eta \biggr\}\cap\operatorname{Ext}^c.
\]
Then
\[
\mathbh{1}_{A_n(K,\eta)}\sum_{\mathbf{i} \in\mathbb{G}_n^*}
Z_{\mathbf{i}}\mathbh{1}_{\{Z_{\mathbf{i}}> K\}}
\geq\mathbh{1}_{A_n(K,\eta)}\mathcal{Z}_n \eta.
\]
Using (\ref{minW}), we have
\[
\mathbh{1}_{A_n(K,\eta)} (2m)^{-n} \sum_{\mathbf{i} \in\mathbb{G}_n^*}
Z_{\mathbf{i}} \mathbh{1}_{\{Z_{\mathbf{i}}> K\}}
\geq\eta\mathbh{1}_{A_n(K,\eta)}C
\]
so that taking expectations,
\begin{eqnarray*}
m^{-n} \mathbb{E}\Biggl(2^{-n} \sum_{\mathbf{i} \in
\mathbb{G}_n^*} Z_{\mathbf{i}} \mathbh{1}_{\{Z_{\mathbf{i}}> K\}}\Biggr)
&\geq& \mathbb{E}\bigl(\mathbh{1}_{A_n(K,\eta)}C\bigr)\eta
\\
m^{-n}\mathbb{E}\bigl(Z_n\mathbh{1}_{\{Z_n>K\}}\bigr)/\eta
&\geq& \mathbb{E}\bigl(\mathbh{1}_{A_n(K,\eta)}C\bigr).
\end{eqnarray*}
Then, by (\ref{control}), we have
\[
\lim_{K\rightarrow\infty} \sup_{n\in\mathbb{N}}
\bigl\{\mathbb{E}\bigl(\mathbh{1}_{A_n(K,\eta)}C\bigr)\bigr\}=0.
\]
Then observe that $\forall\alpha>0$,
$\inf_{\mathbb{P}^*(A)\geq\alpha}\{\mathbb{E}(C\mathbh{1}_{A})\}>0$. So $\exists
K_0\geq 0$ such that $\forall K\geq K_0$, $\forall n \in\mathbb{N}$,
\[
\mathbb{P}^*(A_n(K,\eta))<\alpha,
\]
which completes the proof.
\end{pf}

First, for any $\varepsilon>0$, choose $K$ using the previous lemma such that
\[
\mathbb{P}^*\biggl ( \frac{ \sum_{\mathbf{i} \in\mathbb{G}_n^*}
Z_{\mathbf{i}}\mathbh{1}_{\{Z_{\mathbf{i}}\leq K\}}}{\mathcal{Z}_n} \geq1/2  \biggr)
= 1-\mathbb{P}^* \biggl( \frac{ \sum_{\mathbf{i} \in\mathbb{G}_n^*}
Z_{\mathbf{i}} \mathbh{1}_{\{Z_{\mathbf{i}}> K\}}}{\mathcal{Z}_n}
< 1/2  \biggr)\geq1-\varepsilon/2.
\]
Adding that conditionally on $\operatorname{Ext}^c$, $\mathcal{Z}_n
\stackrel{n\rightarrow\infty}{\longrightarrow} \infty$ a.s, gives
the following result.
\begin{Pte}\label{nbparasup}
Let $\varepsilon>0$, there exists $K \in\mathbb{N}$ such that
$\forall N \in\mathbb{N}$, $\exists n_0 \in\mathbb{N}$
such that $\forall n\geq n_0$,
\[
\mathbb{P}^*\Biggl(\sum_{\mathbf{i} \in\mathbb{G}_n^*} Z_{\mathbf{i}}
\mathbh{1}_{\{Z_{\mathbf{i}}\leq K\}}\geq N\Biggr)\geq1-\varepsilon.
\]
\end{Pte}

Second, we derive an estimation of $\# \mathbb{G}_n^*$. By Lemma \ref{celltc},
the cells are not very contaminated so the number of contaminated
cells is asymptotically proportional to the number of parasites, which
is a Bienaym\'{e} Galton--Watson process.
\begin{Pte}\label{estim}
For every $\varepsilon>0$, there exist $A,B>0$ such that for every $n \in
\mathbb{N}$,
\[
\mathbb{P}^* \biggl( \frac{\# \mathbb{G}_n^*}{(2m)^n} \in[A,B]
\biggr)\geq 1-\varepsilon.
\]
\end{Pte}
\begin{pf}
First use (\ref{minW}) to get
\[
\frac{\# \mathbb{G}_n^*}{(2m)^n} \leq\frac{\mathcal
{Z}_n}{(2m)^n}\leq D.
\]
Moreover using again $(\ref{minW})$, we have
\[
\frac{\# \mathbb{G}_n^*}{(2m)^n}\geq\frac{\sum_{\mathbf{i} \in
\mathbb{G} _n^*} Z_{\mathbf{i}}
\mathbh{1}_{\{Z_{\mathbf{i}}\leq K\}}}{K(2m)^n}\geq\frac{C}{K}\frac
{\sum_{\mathbf{i} \in\mathbb{G}_n^*}
Z_{\mathbf{i}} \mathbh{1}_{\{Z_{\mathbf{i}}\leq K\}}}{\mathcal{Z}_n}
\]
and Lemma \ref{celltc} gives the result.
\end{pf}
%

\subsection{Separation of the descendances of parasites}\label{s6.3}

Start with two parasites and consider the BPRE ($Z_n)_{n\in\mathbb{N}}$.
Even when conditioning on the survival of their descendance,
the descendance of one of them dies out. This ensures
that two distinct parasites in generation $n$ do not have descendants
which belong to the same cell
in generation $n+q$ if $q$ is large enough. More precisely, we define
$N_n(\mathbf{i})$
as the number of parasites of cell $\mathbf{i}\vert n$ whose
descendance is still
alive in cell $\mathbf{i}$ and we prove the following result.
%
\begin{Pte}\label{separation}
$\forall K \in\mathbb{N}$, $\forall \varepsilon$, $\eta>0$, $\exists q\in\mathbb{N}$
such that $\forall n \in\mathbb{N}$, we have
\[
\mathbb{P}^* \biggl( \frac{ \# \{ \mathbf{i} \in\mathbb{G}_{n+q}^*
\dvtx Z_{\mathbf {i}\vert n}\leq K,
N_n(\mathbf{i})\geq2\}}{\# \mathbb{G}_{n+q}^*}\geq\eta \biggr)\leq
\varepsilon.
\]
\end{Pte}
\begin{pf}
Let $K \in\mathbb{N}$, $\eta>0$ and consider for $A>0$,
\[
E_n^q(\eta)= \biggl\{\frac{\# \{ \mathbf{i} \in\mathbb{G}_{n+q}^* \dvtx
Z_{\mathbf{i}\vert n}\leq K,
N_n(\mathbf{i})\geq2\}}{\# \mathbb{G}_{n+q}^* }\geq\eta \biggr\}
\cap \biggl \{ \frac{\# \mathbb{G}_{n+q}^* }{(2m)^{n+q}} \geq A \biggr\}.
\]
Then
\[
\mathbh{1}_{E_n^q(\eta)} \# \{ \mathbf{i} \in\mathbb{G}_{n+q}^* \dvtx
Z_{\mathbf{i}\vert n}\leq K,
N_n(\mathbf{i})\geq2\} \geq\mathbh{1}_{E_n^q(\eta)} \eta A (2m)^{n+q}
\]
so that taking expectations,
\begin{eqnarray*}
\mathbb{P}(E_n^q(\eta))&\leq&
\frac{2^{-(n+q)} \mathbb{E}( \sum_{\mathbf{i} \in\mathbb{G}_{n+q}}
\mathbh{1}_{\{Z_{\mathbf{i}\vert n}\leq K, N_n(\mathbf{i})\geq2 \}}) }
{ \eta A m^{n+q} }
\\
&\leq& \frac{ 2^{-n} \sum_{\mathbf{j} \in\mathbb{G}_n} \mathbb{P}
(0<Z_{\mathbf{j}}\leq K) 2^{-q} \sum_{\mathbf{i} \in\mathbb{G}_q}
\mathbb{P}_K(N_0(\mathbf{i})\geq2) }{\eta A m^{n+q}}
\\
&\leq& \frac{\mathbb{P}(Z_n>0)2^{-q}\sum_{\mathbf{i}
\in\mathbb{G}_q} \mathbb{P}_K
(N_0(\mathbf{i})\geq2) }{\eta A m^{n+q}}.
\end{eqnarray*}
As we have ${K\choose 2}$ ways to choose two parasites among $K$ and
they both survive along $\mathbf{i}$
with probability $\mathbb{P}(Z_{\mathbf{i}}>0)^2$, we have
\[
\mathbb{P}_K\bigl(N_0(\mathbf{i})\geq2\bigr)
\leq \pmatrix{K\cr 2} \mathbb{P}(Z_{\mathbf{i}}>0)^2.
\]
Then
\begin{eqnarray*}
\mathbb{P}(E_n^q(\eta))&\leq& \frac{{K\choose 2}2^{-q}\sum_{\mathbf{i}
\in \mathbb{G}_q} \mathbb{P}(Z_{\mathbf{i}}>0)^2 }{\eta A m^{q}}.
\end{eqnarray*}
Conclude
choosing $A$ in agreement with Proposition $\ref{estim}$ and $q$
with (\ref{negli}).
\end{pf}
%

\subsection{Control of filled-in cells}\label{s6.4}

Here we prove that filled-in cells have asymptotically no impact on the
proportions of cells
with a given number of parasites.
\begin{Lem}\label{contr}
$\forall \varepsilon$, $\eta>0$, $\exists K \in\mathbb{N}$ such that
$\forall n,q \in\mathbb{N}$, we have
\[
\mathbb{P}^* \biggl( \frac{\# \{ \mathbf{i} \in\mathbb{G}_{n+q}^* \dvtx
Z_{\mathbf{i}\vert n}>K\}}{\# \mathbb{G}_{n+q}^* }\geq\eta \biggr)\leq\varepsilon.
\]
\end{Lem}
\begin{pf}
Let $\eta>0$, $A>0$ and consider
\[
F_n^q(\eta)= \biggl\{
\frac{\# \{ \mathbf{i} \in\mathbb{G}_{n+q}^* \dvtx
Z_{\mathbf{i}\vert n}>K\}}{\#\mathbb{G}_{n+q}^*} \geq\eta \biggr\}\cap
\biggl\{ \frac{\#\mathbb{G}_{n+q}^*}{(2m)^{n+q}} \geq A \biggr\}
\]
then
\[
\mathbh{1}_{F_n^q(\eta)} \# \{ \mathbf{i} \in\mathbb{G}_{n+q}^* \dvtx
Z_{\mathbf{i}\vert n}>K\} \geq\mathbh{1}_{F_n^q(\eta)} \eta A (2m)^{n+q}.
\]
Taking expectations leads to
\begin{eqnarray*}
\mathbb{P}(F_n^q(\eta))
&\leq& \frac{2^{-(n+q)} \mathbb{E}( \sum_{\mathbf{i} \in
\mathbb{G}_{n+q}} \mathbh{1}_{\{Z_{\mathbf{i}\vert n}>K, Z_{\mathbf{i}}>0 \}})}
{ \eta A m^{n+q} }
\\
&\leq& \frac{2^{-(n+q)} \sum_{\mathbf{i} \in\mathbb{G}_{n+q}}
\mathbb{P}(Z_{\mathbf{i}\vert n}>K, Z_{\mathbf{i}}>0) }{ \eta A m^{n+q} }
\\
&\leq& \frac{\sum_{k> K} 2^{-n} \sum_{\mathbf{j} \in\mathbb{G}_{n}} \mathbb{P}
(Z_{\mathbf{j}}=k) 2^{-q} \sum_{\mathbf{i} \in\mathbb{G}_q}
\mathbb{P}_k(Z_{\mathbf{i}}>0) }{ \eta A m^{n+q} }.
\end{eqnarray*}
Moreover, $\mathbb{P}_k(Z_{\mathbf{i}}>0)=1-(1-\mathbb{P}(Z_{\mathbf{i}}>0))^k\leq
k\mathbb{P}(Z_{\mathbf{i}}>0)$ and we have
\begin{eqnarray*}
\mathbb{P}(F_n^q(\eta))
&\leq& \frac{\sum_{k> K} 2^{-n} \sum_{\mathbf{j}\in\mathbb{G}_{n}}
k\mathbb{P}(Z_{\mathbf{j}}=k) \mathbb{P}(Z_q>0)}{ \eta A m^{n+q}}
\\
&\leq& \frac{ \mathbb{E}(Z_n \mathbh{1}_{\{Z_n> K\}})}{ \eta Am^{n} }
\qquad\mbox{using (\ref{majP})}.
\end{eqnarray*}
By (\ref{negli}), we get
\[
\lim_{K\rightarrow\infty} \sup_{n\in\mathbb{N}}
\{\mathbb{P}(F_n^q(\eta))\}=0.
\]
Complete the proof choosing $A$ in
agreement with Proposition $\ref{estim}$.
\end{pf}

\subsection{\texorpdfstring{Proof of Theorem \textup{\protect\ref{limf1}}}{Proof of Theorem \textup{5.2}}}\label{s6.5}

Consider the contaminated cells in generation \mbox{$n+q$}. Their ancestors in
generation $n$ are cells which are
not very contaminated (by Lemma \ref{contr}). Then taking $q$ large,
the parasites of a contaminated cell
in generation $n+q$ come from
a same parasite in generation $n$ (separation of the descendances of
parasites, Proposition \ref{separation}). Thus
at generation $n+q$,
everything occurs as if all parasites from generation~$n$
belonged to different cells. As the number of parasites at
generation $n$ tends to infinity ($n\rightarrow\infty, \ m_0+m_1>1$),
we have a law of
large numbers phenomenon and get a deterministic limit.

\textit{Step} 1. We prove that for all $\varepsilon, \eta>0$, there exist $n_0
\in\mathbb{N}$ and $\vec{f} \in\mathbb{S}^1(\mathbb{N})$ such
that for every
$n\geq n_0$,
\[
\mathbb{P}^*\bigl(\| (F_k(n))_{k \in\mathbb{N}} -\vec{f} \|_1\geq\eta\bigr)
\leq\varepsilon.
\]

For every $k \in\mathbb{N}^*$ and every parasite $\mathbf{p}$ in
generation $n$, we denote by $Y_k^{q}(\mathbf{p})$ the
number of cells in generation $n+q$
which contain at least $k$ parasites, exactly~$k$ of which have
$\mathbf{p}$ as an ancestor. By convention, $Y_0^{q}(\mathbf{p})=0$.
That is, writing for $\mathbf{p}$ parasite,
$\mathbf{p}\hookrightarrow\mathbf{i}$ when $\mathbf{p}$ belongs to the
cell $\mathbf{i}$
and $\mathbf{p}\vert n$ its ancestor (parasite) in generation $n$,
\[
Y_k^{q}(\mathbf{p})= \sum_{ \mathbf{i} \in\mathbb{G}_{n+q}}
\mathbh{1}_{ \#\{\mathbf{r} \dvtx  \mathbf{r}\hookrightarrow\mathbf{i},
\mathbf{r}\vert n=\mathbf{p} \}=k}, \qquad
k \in\mathbb{N}^*.
\]
By the branching property, $ (Y_k^{q}(\mathbf{p}))_{k\in\mathbb{N}}$
[$\mathbf{p} \in\mathcal{P}(n)$] are i.i.d. and we denote by
$(Y_k^{q})_{k\in\mathbb{N}}$ a random variable
with this common distribution. Denoting by $\mathcal{P}_{K}(n)$ the
set of parasites in generation
$n$ which belong to a cell containing at most $K$ parasites, we have
%
\begin{eqnarray}\label{ineg}
&&\quad \sum_{k \in\mathbb{N}^*}  \Biggl\vert\# \{ \mathbf{i}\in\mathbb{G}
_{n+q}^* \dvtx Z_{\mathbf{i}}=k\}
- \sum_{\mathbf{p}\in\mathcal{P}_{K}(n)} Y^{q}_k(\mathbf{p})  \Biggr\vert
\nonumber\\[-8pt]
\\[-8pt]
&&\quad\qquad \leq (K+1) \# \{ \mathbf{i} \in\mathbb{G}_{n+q}^* \dvtx
Z_{\mathbf{i}\vert n}\leq K,  N_n(\mathbf{i})\geq2\}
+ \# \{ \mathbf{i} \in\mathbb{G}_{n+q}^* \dvtx Z_{\mathbf{i}\vert n}>K\}.
\nonumber\hspace*{-18pt}
\end{eqnarray}
Indeed, the left-hand side of (\ref{ineg}) is less than
\[
\sum_{k \in\mathbb{N}^*}  \Biggl\vert\# \{ \mathbf{i}\in\mathbb{G}_{n+q}^* \dvtx
Z_{\mathbf{i}}=k, Z_{\mathbf{i}\vert n}\leq K \}
- \sum_{\mathbf{p}\in\mathcal{P}_{ K}(n)} Y^{q}_k(\mathbf{p})\Biggr \vert
+ \# \{\mathbf{i} \in\mathbb{G}_{n+q}^*\dvtx
Z_{\mathbf{i}\vert n}>K\}.
\]
And recalling that $N_n(\mathbf{i})$ is the number of parasites of
cell $\mathbf{i}\vert n$ whose descendance is still
alive in cell $\mathbf{i}$, we get the following equalities:
\begin{eqnarray*}
\sum_{\mathbf{p}\in\mathcal{P}_{ K}(n)} Y^{q}_k(\mathbf{p})
&=& \sum_{ \mathbf{i} \in\mathbb{G}_{n+q}}
\sum_{\mathbf{p} \in\mathcal{P}_{ K}(n)} \mathbh{1}_{ \#\{\mathbf{r} \dvtx
\mathbf{r}\hookrightarrow\mathbf{i}, \mathbf{r}\vert n=\mathbf{p} \}=k}
\end{eqnarray*}
and
\begin{eqnarray*}
\mathbh{1}_{Z_{\mathbf{i}}=k,  Z_{\mathbf{i}\vert n}\leq K,
N_n(\mathbf{i})= 1}&=& \mathbh{1}_{N_n(\mathbf{i})= 1}
\sum_{\mathbf{p} \in\mathcal{P}_{ K}(n)} \mathbh{1}_{ \#\{\mathbf{r} \dvtx
\mathbf{r}\hookrightarrow\mathbf{i}, \mathbf{r}\vert n=\mathbf{p} \}=k}
\end{eqnarray*}
which ensure
\begin{eqnarray*}
&& \sum_{k \in\mathbb{N}^*}  \Biggl\vert\# \{ \mathbf{i}\in\mathbb{G}_{n+q}^* \dvtx
Z_{\mathbf{i}}=k,  Z_{\mathbf{i}\vert n}\leq K \}
- \sum_{\mathbf{p}\in\mathcal{P}_{ K}(n)} Y^{q}_k(\mathbf{p})\Biggr\vert
\\
&&\qquad \leq \sum_{k \in\mathbb{N}^*} \sum_{{\mathbf{i} \in\mathbb{G}
_{n+q}, N_n(\mathbf{i})\geq2}}
\Biggl\vert\mathbh{1}_{Z_{\mathbf{i}}=k,
Z_{\mathbf{i}\vert n}\leq K} - \sum_{\mathbf{p} \in\mathcal{P}_{
K}(n)}\mathbh{1}_{ \#\{\mathbf{r} \dvtx  \mathbf{r}\hookrightarrow\mathbf{i},
\mathbf{r}\vert n=\mathbf{p} \}= k}\Biggr  \vert
\\
&&\qquad \leq \# \{ \mathbf{i} \in\mathbb{G}_{n+q}^* \dvtx
Z_{\mathbf{i}\vert n}\leq K, \ N_n(\mathbf{i})\geq2\}
+ \mathop{\sum_{\mathbf{i} \in\mathbb{G}_{n+q}, N_n(\mathbf{i})\geq2}}
_{\mathbf{p} \in\mathcal{P}_{ K}(n)} \mathbh{1}_{ \#\{\mathbf{r}
\dvtx  \mathbf{r}\hookrightarrow\mathbf{i},
\mathbf{r}\vert n=\mathbf{p} \}>0}
\\
&&\qquad \leq \# \{ \mathbf{i} \in\mathbb{G}_{n+q}^* \dvtx
Z_{\mathbf{i}\vert n}\leq K, \ N_n(\mathbf{i})\geq2\}
+ \sum_{{\mathbf{i} \in\mathbb{G}_{n+q}, N_n(\mathbf{i})\geq 2}}
K \mathbh{1}_{Z_{\mathbf{i}\vert n}\leq K}
\\
&&\qquad = (K+1)\# \{ \mathbf{i} \in\mathbb{G}_{n+q}^* \dvtx
Z_{\mathbf{i}\vert n}\leq K,  N_n(\mathbf{i})\geq2\}.
\end{eqnarray*}

We shall now prove that the quantities on the right-hand side of
(\ref{ineg}) are small when $n$ and $q$
are large enough and that
$\sum_{\mathbf{p}\in\mathcal{P}_{ K}(n)} Y^{q}_k(\mathbf{p})$
follow a law of large number.
To that purpose, let $\varepsilon,\eta>0$ and for all $K,k,n,q\geq0$ define
%
\[
G^K_k(n,q):= \frac{ \sum_{\mathbf{p}\in\mathcal{P}_{ K}(n)}
Y^{q}_k(\mathbf{p}) }{\sum_{k \in\mathbb{N}}
\sum_{\mathbf{p}\in\mathcal{P}_{ K}(n)} Y^{q}_k(\mathbf{p})}.
\]


First, by Proposition \ref{nbparasup} and (\ref{P_n}),
$\exists K_1 \in\mathbb{N}$
such that $\forall N \in\mathbb{N}$, $\exists n_1 \in\mathbb{N}$
such that $\forall
K\geq K_1$, $\forall n\geq n_1$,
%
\begin{equation}\label{lgn}
\mathbb{P}^n\bigl(\vert\mathcal{P}_{ K}(n)\vert\geq N\bigr) \geq1-\varepsilon.
\end{equation}
Moreover by Lemma $\ref{contr}$, $\exists  K_2\geq K_1$ such that
$\forall n,q \in\mathbb{N}$,
%
\begin{equation}\label{ineg1}
\mathbb{P}^* \biggl(\frac{ \# \{ \mathbf{i} \in\mathbb{G}_{n+q}^* \dvtx
Z_{\mathbf{i}\vert n}>K_2\}}{\# \mathbb{G}_{n+q}^* }\geq\eta \biggr)
\leq \varepsilon.
\end{equation}
And by Proposition \ref{separation}, $\exists q_0 \in\mathbb{N}$
such that $\forall n \in\mathbb{N}$,
%
\begin{equation}\label{ineg2}
\mathbb{P}^*\biggl(\frac{ \# \{ \mathbf{i} \in\mathbb{G}_{n+q_0}^* \dvtx
Z_{\mathbf{i}\vert n}\leq K_2,
N_n(\mathbf{i})\geq2\}}{\# \mathbb{G}_{n+q_0}^* }\geq\eta/(K_2+1)\biggr)
\leq\varepsilon.
\end{equation}
Use then (\ref{ineg}), (\ref{ineg1}) and (\ref{ineg2}) to get
\[
\mathbb{P}^* \biggl(\frac{\sum_{k \in\mathbb{N}^*}  \vert\# \{
\mathbf{i}\in \mathbb{G}_{n+q_0}\dvtx
Z_{\mathbf{i}}=k\} - \sum_{\mathbf{p}\in\mathcal{P}_{K_2}(n)}
Y^{q_0}_k(\mathbf{p})\vert}{\# \mathbb{G}^*_{n+q_0}}
\geq2\eta\biggr )\leq2\varepsilon.
\]
Then by (\ref{impliquer}), for every $n \in\mathbb{N}$, we have
%
\begin{equation}\label{cvv}
\mathbb{P}^* \bigl( \bigl\| \bigl(F_k(n+q_0)\bigr)_{k\in\mathbb{N}}
- \bigl(G^{K_2}_k(n,q_0) \bigr)_{k \in\mathbb{N}} \bigr\|_1\geq4\eta
\bigr)\leq2\varepsilon.
\end{equation}

Second, conditionally on $\mathcal{Z}_n>0$,
$Y_k^{q_0}(\mathbf{p})$ [$\mathbf{p} \in\mathcal{P}_{K_2}(n)$]
are i.i.d. Then the law of large numbers (LLN)
ensures that $\forall k \in
\mathbb{N}$, as $n$ and so $\mathcal{P}_{K_2}(n)$ becomes large:
\[
G^{K_2}_k(n,q_0) \longrightarrow f_k(q_0) \qquad\mbox{where }
f_k(q_0):= \frac{ \mathbb{E}(Y_k^{q_0}) }
{\sum_{k' \in\mathbb{N}} \mathbb{E}(Y_{k'}^{q_0})}.
\]
To see that, divide the numerator and denominator of $G^{K_2}_k(n,q_0)$
by $\#\mathcal{P}_{K_2}(n)$. More precisely, by the LLN,
there exists $N>0$ such that for all $n \in\mathbb{N}$,
\[
\mathbb{P}^n \bigl( \| (G^{K_2}_k(n,q_0) )_{k \in\mathbb{N}^*}
- \vec{f}(q_0) \|_1\geq\eta,
\mathcal{P}_{K_2}(n)\geq N  \bigr)\leq\varepsilon.
\]
%
So using (\ref{lgn}), there exists $n_1 \in\mathbb{N}$ such that
for every $\forall n\geq n_1$,
\[
\mathbb{P}^n
\bigl( \| (G^{K_2}_k(n,q_0) )_{k \in\mathbb{N}^*}
- \vec{f}(q_0) \|_1\geq\eta \bigr)\leq2\varepsilon.
\]
Finally by $(\ref{P_n})$, there exists $n_2\geq n_1$ such that for
every $n\geq n_2$,
%
\begin{equation}\label{LBN}
\mathbb{P}^* \bigl( \| (G^{K_2}_k(n,q_0) )_{k \in\mathbb{N}^*}
- \vec{f}(q_0) \|_1\geq\eta \bigr)\leq3\varepsilon.
\end{equation}

As a conclusion, using $(\ref{cvv})$ and $(\ref{LBN})$, we have
proved that for all $\varepsilon, \eta>0$,
and for every $n\geq n_2+q_0$,
\[
\mathbb{P}^*\bigl(\| (F_k(n))_{k \in\mathbb{N}^*}
- \vec{f}(q_0) \|_1\geq 5\eta\bigr)\leq 3\varepsilon.
\]

\textit{Step} 2. Existence of the limit.

For every $l \in\mathbb{N}$, there exist $n_0(l) \in\mathbb{N}$ and
$\vec{f} (l)
\in$ $\mathbb{S}^1(\mathbb{N})$ such that for every $n\geq n_0(l)$
\[
\mathbb{P}\bigl( \| F(n)-\vec{f}(l) \|_1\geq1/2^{l+1}\bigr)
\leq1/2^{l}.
\]
Then for all $l,l'$ such that $2\leq l\leq l'\dvtx \|\vec{f}(l')
-\vec{f}(l) \|_1\leq1/2^{l}$ and completeness of $l^1(\mathbb{N})$ ensures
that $(\vec{f} (l))_{l\in\mathbb{N}}$ converges in
$\mathbb{S}^1(\mathbb{N})$ to a
limit $\vec{f}$. Moreover, $\|\vec{f}(l) -\vec{f} \|_1\leq1/2^{l}$
so for every $n\geq n_0(l)$,
\[
\mathbb{P}\bigl( \| F(n)- \vec{f} \|_1\geq1/2^{l}\bigr)\leq1/2^{l}
\]
which ensures the convergence in probability of $( F_k(n))_{n\in\mathbb{N}}$
to $\vec{f}$ as $n \rightarrow\infty$.

\textit{Step} 3.
Characterization of the limit as $f_k=\mathbb{P}(\Upsilon=k)$.

By Proposition $\ref{yag}$, we have
%
\begin{equation}\label{limbp}
\forall k \in\mathbb{N}\qquad
\mathbb{P}(Z_n=k  \mid  Z_n \ne0)\stackrel{n\rightarrow\infty}{\longrightarrow}
\mathbb{P}(\Upsilon=k).
\end{equation}

Moreover, for every $k \in\mathbb{N}^*$, using (\ref{kim}),
\[
\mathbb{P}(Z_n=k  \mid  Z_n \ne0)
= \frac{\mathbb{E}( \# \{\mathbf{i} \in\mathbb{G}_n \dvtx
Z_{\mathbf{i}} =k \})}{\mathbb{E}(\# \mathbb{G}^*_n )}
= \frac{\mathbb{E}( F_k(n) \#\mathbb{G}_n^*)}{\mathbb{E}(\# \mathbb{G}^*_n)}.
\]
As $F_k(n)$ converges in probability
to a deterministic limit $f_k$, we get
%
\begin{equation}\label{fre}
\forall k \in\mathbb{N} \qquad
\mathbb{P}(Z_n=k  \mid  Z_n \ne0)
\stackrel{n\rightarrow\infty}{\longrightarrow} f_k.
\end{equation}
Indeed, by
Proposition \ref{estim}, there exists $A>0$ such that
\[
\frac{\mathbb{E}(\# \mathbb{G}_n^*)}{(2m)^n}\geq A.
\]
Then for every $\eta>0$, using $| F_k(n)-f_k|\leq1$, we have
\begin{eqnarray*}
\biggl \vert\frac{\mathbb{E}( F_k(n) \# \mathbb{G}_n^*)}
{\mathbb{E}(\# \mathbb{G}^*_n )} -f_k\biggr\vert
&\leq& \frac{\mathbb{E}( \# \mathbb{G}_n^* \mid F_k(n)-f_k\mid\mathbh{1}_{\{\mid
F_k(n)-f_k\mid< \eta\}})}{\mathbb{E}( \# \mathbb{G}_n^* )}
\\
&&{} +  \frac{\mathbb{E}( \# \mathbb{G}_n^* \mathbh{1}_{\{| F_k(n)-f_k|\geq\eta\}})}
{\mathbb{E}( \# \mathbb{G}_n^* )}
\\
&\leq& \eta+\frac{\mathbb{E}(\mathcal{Z}_n \mathbh{1}_{\{|F_k(n)-f_k|\geq\eta\}})}
{A(2m)^n}.
\end{eqnarray*}
By (\ref{mom2GW}), $\mathcal{Z}_n/(2m)^n$ is
bounded in $L^2$ and it is uniformly integrable.
Then, thanks to the previous steps,
the second term in the last displayed equation vanishes as $n$ grows
and we get (\ref{fre}). Putting
($\ref{limbp}$) and ($\ref{fre}$) together proves that $f_k=\mathbb{P}
(\Upsilon=k)$.

\subsection{Proof of corollaries}\label{s6.6}
\mbox{}
\begin{pf*}{Proof of Corollary \ref{gn}}
Recall that $\mathbb{E}(\Upsilon)<\infty$ (Proposition \ref{yag})
and note also that for every $K \in\mathbb{N}^*$,
\[
\# \mathbb{G}_n^*=\frac{\sum_{\mathbf{i} \in\mathbb{G}_n^*}
Z_{\mathbf{i}}
\mathbh{1}_{\{Z_{\mathbf{i}}\leq K\}}}{ \sum_{k=1}^K kF_k(n)}.
\]
Then using
$\sum_{\mathbf{i} \in\mathbb{G}_n^*} Z_{\mathbf{i}}
\mathbh{1}_{\{Z_{\mathbf{i}}\leq K\}} \leq\mathcal{Z}_n$ gives
\begin{eqnarray*}
\biggl \vert\frac{\# \mathbb{G}_n^*}{\mathcal{Z}_n}
- \frac{1}{\mathbb{E}(\Upsilon )} \biggr\vert
&=&  \biggl\vert\frac{1}{\sum_{k=1}^K kF_k(n)}
\frac{\sum_{\mathbf{i}\in\mathbb{G}_n^*} Z_{\mathbf{i}}
\mathbh{1}_{\{Z_{\mathbf{i}}\leq K\}}}{\mathcal{Z}_n}
- \frac{1}{\mathbb{E}(\Upsilon )}  \biggr\vert
\\
&\leq&  \biggl\vert\frac{1}{\sum_{k=1}^K kF_k(n)}
- \frac{1}{\mathbb{E}(\Upsilon)}  \biggr\vert
+ \frac{1}{\mathbb{E}(\Upsilon)}
\biggl\vert\frac{\sum_{\mathbf{i} \in\mathbb{G}_n^*} Z_{\mathbf{i}}
\mathbh{1}_{\{Z_{\mathbf{i}}\leq K\}}}{\mathcal{Z}_n}-1 \biggr\vert.
\end{eqnarray*}
Let $\eta,\varepsilon>0$. We use Lemma \ref{celltc} to choose $K \in
\mathbb{N}^*$ such that
\begin{eqnarray*}
\forall n \in\mathbb{N} \qquad
\mathbb{P}^* \biggl( \frac{ \sum_{\mathbf{i} \in\mathbb{G}_n^*}
Z_{\mathbf{i}} \mathbh{1}_{\{Z_{\mathbf{i}}\leq K\}}}
{\mathcal{Z}_n}\geq1-\eta \biggr)& \geq& 1-\varepsilon; \qquad
\\
\biggl\vert\frac{1}{ \mathbb{E}(\Upsilon\mathbh{1}_{\Upsilon\leq K}) }
-\frac{1}{\mathbb{E}(\Upsilon )} \biggr \vert &\leq& \eta.
\end{eqnarray*}
Choose $n_0 \in\mathbb{N}$ using Theorem \ref{limf1} so that for every
$n\geq n_0$,
\[
\mathbb{P}^* \biggl( \biggl\vert\frac{1}{\sum_{k=1}^K kF_k(n)}
- \frac{1}{\mathbb{E}(\Upsilon \mathbh{1}_{\Upsilon \leq K})}
\biggr\vert\leq\eta \biggr)\geq1-\varepsilon.
\]
Then for every $n\geq n_0$,
\[
\mathbb{P}^* \biggl( \biggl\vert\frac{ \# \mathbb{G}_n^*}{\mathcal
{Z}_n}-\frac{1}{\mathbb{E}(\Upsilon
)} \biggr\vert\geq2\eta+\frac{1}{\mathbb{E}(\Upsilon )}\eta
\biggr)\leq2\varepsilon,
\]
which proves the convergence in probability of
$\#\mathbb{G}_n^*/\mathcal{Z}_n$ to
$1/\mathbb{E}(\Upsilon )$. The second convergence follows from
(\ref{nbparas}).
\end{pf*}
\begin{pf*}{Proof of Corollary \ref{limf13}}
We write for $n,p,k$ $\in\mathbb{N}$,
\[
\frac{\#\{\mathbf{i} \in\mathbb{G}_{n+p}^* \dvtx
Z_{\mathbf{i}\mid n}=k\}}{\#\mathbb{G}^*_{n+p}}
= \frac{(2m)^p}{\#\mathbb{G}^*_{n+p}} \sum_{\mathbf{j}
\in\mathbb{G}_n^* \dvtx Z_{\mathbf{j}}=k}
\frac{\#\{\mathbf{i}\in\mathbb{G}_{n+p}^* \dvtx
\mathbf{i}\vert n=\mathbf{j}\}}{(2m)^p}.
\]
Conditionally on $Z_{\mathbf{j}}=k$, by Corollary \ref{gn} and
separation of descendances of parasites, we have the following
convergence in probability
\[
\frac{ \#\{\mathbf{i}\in\mathbb{G}_{n+p}^* \dvtx
\mathbf{i}\vert n=\mathbf{j}\}}{(2m)^p}
\stackrel{p\rightarrow\infty}{\longrightarrow}W_k(\mathbf{j}),
\]
where $W_k(\mathbf{j})$ is the sum of $k$ i.i.d. variables distributed
as $W/\mathbb{E}(\Upsilon )$. Then, using also (\ref{nbparas}),
%
\begin{equation}\label{lln}
\mathbb{E}(W_k(\mathbf{j}))
= \frac{k\mathbb{E}(W)}{\mathbb{E}(\Upsilon )}
= \frac{k}{\mathbb{E}(\Upsilon )}.
\end{equation}
Using again Corollary \ref{gn}, we get the first limit of the corollary
\[
\lim_{p\rightarrow\infty} \frac{\#\{\mathbf{i} \in\mathbb{G}_{n+p}^* \dvtx
Z_{\mathbf{i}\mid n}=k\}}{\#\mathbb{G}^*_{n+p}}
\stackrel{\mathbb{P}}{=}
\frac{\mathbb{E}(\Upsilon )}{W} \frac{\sum_{\mathbf{j} \in \mathbb{G}_n^* \dvtx
Z_{\mathbf{j}}=k}W_k(\mathbf{j})}{(2m)^n}.
\]

Moreover, Theorem \ref{limf1} ensures that
\[
\frac{\#\{\mathbf{j} \in\mathbb{G}_n^* \dvtx Z_{\mathbf{j}}=k\}}
{(2m)^n}=F_k(n)\frac{\mathcal{Z}_n}{(2m)^n} \stackrel{n\rightarrow
\infty}{\longrightarrow} \frac{W}{\mathbb{E}(\Upsilon )} f_k.
\]
And conditionally on $\#\mathbb{G}_n^*>0$, $ W_k(\mathbf{j})$
($\mathbf{j}\in\mathbb{G}_n^*$) is i.i.d. by the branching property and $\#
\mathbb{G}_n^*$
tends to infinity.
So the law of large numbers and (\ref{lln}) ensure that
\begin{eqnarray*}
&& \lim_{n\rightarrow\infty}\frac{\mathbb{E}(\Upsilon )}{W}
\frac{\sum_{\mathbf{j} \in\mathbb{G}_n^* \dvtx
Z_{\mathbf{j}}=k}W_k(\mathbf{j})}{(2m)^n}
\\
&&\qquad = \lim_{n\rightarrow\infty}\frac{\mathbb{E}(\Upsilon )}{W}
\frac{\#\{\mathbf{j} \in\mathbb{G}_n^* \dvtx Z_{\mathbf{j}}=k\}
}{(2m)^n}\frac{\sum_{\mathbf{j} \in\mathbb{G}_n^* \dvtx
Z_{\mathbf{j}}=k}W_k(\mathbf{j})}{\#\{\mathbf{j} \in\mathbb{G}_n^* \dvtx
Z_{\mathbf{j}}=k\}}
\stackrel{\mathbb{P}^*}{=} \frac{kf_k}{\mathbb{E}(\Upsilon )},
\end{eqnarray*}
which ends the proof.
\end{pf*}

\section{Proofs in the case $(m_0,m_1)\in D_1$}\label{s7}

We still assume $\mathbb{E} (Z^{(a)2} )<\infty$, the proof is
in the same vein as the proof in the previous section and use the
separation of the descendances of the parasites. The main
difference
is that $\mathcal{Z}_n$ does not explode so the limit is not
deterministic and
the convergence holds in distribution.
\begin{Lem}\label{separationd}
For every $K>0$, there exists $q_0\in\mathbb{N}$ such that for all
$q\geq q_0$ and $n\in\mathbb{N}$,
\[
\mathbb{P}^{n+q} \bigl( \{ \mathbf{i} \in\mathbb{G}_{n+q}^* \dvtx
N_n(\mathbf{i})\geq2\}\ne \varnothing,
\mathcal{Z}_{n}\leq K  \bigr)\leq\varepsilon.
\]
\end{Lem}
\begin{pf}
Denoting by $E_n^q$ the event
\[
\bigl\{\{ \mathbf{i} \in\mathbb{G}_{n+q}^* \dvtx
N_n(\mathbf{i})\geq2\} \ne \varnothing,
\mathcal{Z}_{n}\leq K\bigr\},
\]
we have
\[
\mathbh{1}_{E_n^q}\leq\sum_{\mathbf{i}\in\mathbb{G}_{n+q}}
\mathbh{1}_{\{N_n(\mathbf{i})\geq2,
\mathcal{Z}_{n}\leq K\}}.
\]
Thus we can follow the proof of Lemma $\ref{separation}$.
\begin{eqnarray*}
\mathbb{P}^{n+q}(E_n^q )
&\leq& \sum_{\mathbf{i}\in\mathbb{G}_{n+q}}
\frac{ \mathbb{P}(N_n(\mathbf{i})\geq2,
\mathcal{Z}_{n}\leq K)}{\mathbb{P}(\mathcal{Z}_{n+q}>0)}
\\
&\leq& \frac{\sum_{\mathbf{i}\in\mathbb{G}_{n+q}} \mathbb{P}
(N_n(\mathbf{i})\geq2 ,
Z_{\mathbf{i}\mid n}\leq K)} {U(2m)^{n+q}} \qquad\mbox{using (\ref{mintot})}
\\
&\leq& \frac{ \mathbb{P}(0<Z_n\leq K) 2^{-q}\sum_{i\in\mathbb{G}_q}
\mathbb{P}_K(N_0(\mathbf{i})\geq2)}{U m^{n+q}}
\\
&\leq& \frac{{K\choose 2}2^{-q}\sum_{i\in\mathbb{G}_q}
\mathbb{P}(Z_{\mathbf{i}}>0)^2}{Um^{q}}
\qquad\mbox{using (\ref{majP}).}
\end{eqnarray*}
Conclude with (\ref{negli}).
\end{pf}
\begin{pf*}{Proof of Theorem \ref{limf2}}
\mbox{}

\textit{Step} 1. We recall that $\mathcal{P}_n$ is the set of parasites in
generation $n$, follow Step 1 in the proof
of Theorem \ref{limf1} and use its notation. Thus, we begin with
proving that for every $\varepsilon>0$, there exists
$n_0 \in\mathbb{N}$ such that for every $n\geq n_0$,
\[
\mathbb{P}^{n+q} \bigl(\bigl\|\bigl(\#\{\mathbf{i} \in\mathbb{G}_{n+q}^* \dvtx
Z_{\mathbf{i}}=k \}\bigr)_{k \in\mathbb{N}}
- (N_k(n,q))_{k\in\mathbb{N}} \bigr\|_1 \ne 0  \bigr)\leq\varepsilon,
\]
where for all
$n,q,k\geq0$, $N_k(n,q):= \sum_{\mathbf{p}\in\mathcal{P}(n)} Y^{q}_k(\mathbf{p})$.

First, by (\ref{limlim}), there exist $K,q_0$ $\in\mathbb{N}$ such that
for every $q\geq q_0$,
%
\begin{equation}\label{limnq}
\lim_{n\rightarrow\infty} \mathbb{P}^{n+q}(\mathcal{Z}_n>K)\leq\varepsilon.
\end{equation}
By Lemma \ref{separationd}, there exists $q_1\geq q_0$ such that for every
$n \in\mathbb{N}$, we have
%
\begin{equation}\label{separadis}
\mathbb{P}^{n+q_1} \bigl( \{\mathbf{i} \in\mathbb{G}_{n+q_1}^* \dvtx
N_n(\mathbf{i})\geq2\} \ne \varnothing,
\mathcal{Z}_{n}\leq K  \bigr)\leq\varepsilon.
\end{equation}
And by $(\ref{limnq})$, there
exists $n_0\geq0$ such that for every $n\geq n_0$,
\[
\mathbb{P}^{n+q_1}(\mathcal{Z}_n\geq K)\leq2\varepsilon.
\]
Then
\[
\mathbb{P}^{n+q_1}
\bigl( \# \{ \mathbf{i} \in\mathbb{G}_{n+q_1}^* \dvtx
N_n(\mathbf{i})\geq2\} \ne0  \bigr)\leq3\varepsilon.
\]
Moreover,
\begin{eqnarray*}
&& \# \{ \mathbf{i} \in\mathbb{G}_{n+q_1}^* \dvtx
N_n(\mathbf{i})\geq2\} =0
\\
&&\quad \quad \Longrightarrow \quad
\bigl(\#\{\mathbf{i} \in\mathbb{G}_{n+q_1}^* \dvtx
Z_{\mathbf{i}}=k\}\bigr)_{k \in\mathbb{N}}=(N_k(n,q_1))_{k\in
\mathbb{N}}.
\end{eqnarray*}
%
Then for every $n\geq n_0$,
\[
\mathbb{P}^{n+q_1} \bigl(\|(\#\{\mathbf{i} \in\mathbb{G}_{n+q_1}^* \dvtx
Z_{\mathbf{i}}=k \})_{k \in\mathbb{N}}-(N_k(n,q_1))_{k\in \mathbb{N}}
\|_1 \ne0 \bigr)\leq3 \varepsilon.
\]

\textit{Step} 2. As $l^1(\mathbb{N})$ is separable, we can consider the distance
$d$ associated with the weak convergence of probabilities on
$l^1(\mathbb{N})$. It is defined for any $\mathbb{P}_1$ and
$\mathbb{P}_2$ probabilities by
(see Theorem 6.2, Chapter II in \cite{par})
\[
d(\mathbb{P}_1,\mathbb{P}_2)=\sup
\biggl\{ \biggl\vert\int f(w) \mathbb{P}_1(dw)
- \int f(w)\mathbb{P}_2(dw)  \vert\dvtx \| f\|_{\infty}\leq1,
\| f\|_{\mathrm{Lips}}\leq 1 \biggr\}
\]
where
\[
\| f\| _{\mathrm{Lips}}=\sup
 \biggl\{\frac{\| f(x)-f(y)\| _1}{\| x-y\|_1}\dvtx
x,y \in\mathbb{S}^1(\mathbb{N}), x\ne y \biggr\}.
\]
We prove now that for every $l\geq1$, there exist $n_0(l)\in\mathbb{N}$
and a measure $\mu(l)$ on $\mathbb{N}^*$ such that for every $n\geq
n_0(l)$,
%
\begin{equation}\label{dcauchy}
d \bigl(\mathbb{P}^n \bigl((\#\{\mathbf{i}\in\mathbb{G}_{n}^* \dvtx
Z_{\mathbf{i}}=k \})_{k \in\mathbb{N}} \in \cdot \bigr),
\mu(l) \bigr)\leq1/2^l.
\end{equation}

For that purpose, let $l \in\mathbb{N}$. By Step~1, choose $q,n_0 \
\in\mathbb{N}$
such that
%
\begin{eqnarray}\label{compdone}
&& \forall n\geq n_0\qquad
d \bigl(\mathbb{P}^{n+q}
\bigl( (\#\{\mathbf{i} \in\mathbb{G}_{n+q}^*\dvtx
Z_{\mathbf{i}}=k \})_{k\in\mathbb{N}} \in\cdot  \bigr),
\nonumber\\[-8pt]
\\[-8pt]
&&\hspace*{123pt}
\mathbb{P}^{n+q}  \bigl((N_k(n,q))_{k\in\mathbb{N}} \in\cdot
 \bigr)  \bigr)\leq1/2^{l+1}.
 \nonumber
\end{eqnarray}
Recall that
$(Y^{q}_k(\mathbf{p}))_{k\in\mathbb{N}} (\mathbf{p}\in\mathcal{P}(n))$
is an i.i.d. sequence distributed as $(Y_k^{q})_{k\in\mathbb{N}}$
and $\# \mathcal{P}(n)=\mathcal{Z}_n$. Thus, under
$\mathbb{P}^{n+q}$, $N_k(n,q)$
is the sum of $\mathcal{Z}_{n}$ variables which are i.i.d., distributed as
$Y_k^{q}$ and independent of $\mathcal{Z}_n$, conditionally
on $\sum_{k \in\mathbb{N}}\sum_{\mathbf{p} \in\mathcal{P}(n)}Y^{q}_k(\mathbf{p})>0$.

Moreover $\mathbb{P}^{n+q}(\mathcal{Z}_n \in\cdot )$
converges weakly as $n \rightarrow\infty$ to a probability $\nu$
[see (\ref{limunq})] and we denote
by $\mathcal{N}$ a random variable with distribution $\nu$ and by
$(Y^{q}_k(i))_{k\in\mathbb{N}} (i\in\mathbb{N})$
an i.i.d. sequence independent of
$\mathcal{N}$ and distributed as $(Y_k^{q})_{k\in\mathbb{N}}$. Then
we have
for $n$ large enough,
%
\begin{equation}\label{compdtwo}
d \bigl(\mathbb{P}^{n+q} \bigl((N_k(n,q))_{k\in\ \mathbb{N}}
\in\cdot \bigr), \mu(l) \bigr)\leq1/2^l,
\end{equation}
where $\mu(l)$ is the
distribution of $ (\sum_{1\leq i \leq\mathcal{N}}
Y^{q}_k(i)  )_{k\in \mathbb{N}}$ conditionally on\break
$\sum_{k\in\mathbb{N}}\sum_{1\leq i\leq\mathcal{N}}Y^{q}_k(i)>0$.
Combining (\ref{compdone}) and (\ref{compdtwo}) gives (\ref{dcauchy}).

\textit{Conclusion}.
As $l^1(\mathbb{N})$ is complete, the space of probabilities
on $l^1(\mathbb{N})$ endowed with $d$ is complete (see Theorem 6.5,
Chapter II in \cite{par}), $(\mu(l))_{l \in\mathbb{N}}$ converges and we get the
convergence of Theorem \ref{limf2}.

We now prove that $\mathbb{E}(\sum_{k\in\mathbb{N}^*}kN_{k})<\infty$. For all
$n,K>0$, we have
\[
\mathbb{E}\Biggl(\sum_{k\geq K} k\#\{ \mathbf{i} \in\mathbb{G}_n^* \dvtx
Z_{\mathbf{i}}=k\} \vert\mathcal{Z}_n>0\Biggr)
\leq\mathbb{E}\bigl(\mathcal{Z}_n \mathbh{1}_{\{\mathcal{Z}_n\geq K\}}
\vert\mathcal{Z}_n>0\bigr)
\leq\frac{\mathbb{E}(\mathcal{Z}_n^2)}{\mathbb{P}(\mathcal{Z}_n>0)K}
\]
which converges uniformly to $0$ as $K\rightarrow\infty$ using
(\ref{mom2GW}). Moreover, Theorem \ref{limf2} and
$k\#\{ \mathbf{i} \in\mathbb{G}_n^* \dvtx Z_{\mathbf{i}}=k\}\leq\mathcal{Z}_n$
ensure that
\[
\lim_{n\rightarrow\infty} \mathbb{E}\Biggl(\sum_{1\leq k\leq K} k\#\{
\mathbf{i} \in\mathbb{G}_n^*\dvtx
Z_{\mathbf{i}}=k\} \mid
\mathcal{Z}_n>0\Biggr)
= \mathbb{E}\Biggl(\sum_{1\leq k\leq K} kN_k\Biggr).
\]
Thus we get the expected limit
\[
\mathbb{E}\Biggl(\sum_{k\in\mathbb{N}} k \#\{ \mathbf{i} \in\mathbb{G}_n^*
\dvtx Z_{\mathbf{i}}=k\} \mid \mathcal{Z}_n>0\Biggr)
\stackrel{n\rightarrow\infty}{\longrightarrow}
\mathbb{E}\Biggl(\sum_{k\in\mathbb{N}^*} kN_k\Biggr)
\]
and recalling Section 2.1, we have also
\[
E\Biggl(\sum_{k\in\mathbb{N}^*} k \#\{ \mathbf{i} \in\mathbb{G}_n^* \dvtx
Z_{\mathbf{i}}=k\} \mid \mathcal{Z}_n>0\Biggr)
=\mathbb{E}(\mathcal{Z}_n \mid \mathcal{Z}_n>0)
\stackrel{n\rightarrow\infty}{\longrightarrow} \mathcal{B}'(1)<\infty.
\]
This completes the proof.
\end{pf*}

The proofs of the corollaries follow those of the
previous section.

\section*{Acknowledgment}

I am very grateful to Amaury Lambert who introduced me
to this topic. This work has largely benefited
from his pedagogical qualities and various suggestions.

\printaddresses

\end{document}